\documentclass[a4paper,10pt]{article}

\usepackage[latin1]{inputenc}
\usepackage{amsmath}
\usepackage{amsfonts}
\usepackage{amssymb}
\usepackage{amsthm}
\usepackage[english]{babel}
\allowdisplaybreaks[2]
\usepackage{paralist}
\usepackage[colorlinks=true]{hyperref}
\hypersetup{urlcolor=blue, citecolor=red}
\usepackage{hyperref}

\author{Bj\"orn Augner\footnote{School of Mathematics and Natural Sciences, University of Wuppertal, Gau\ss{}stra\ss{}e 20, D-42119 Wuppertal, Germany. Support by Deutsche Forschungsgemeinschaft (Grant JA 735/8-1) is gratefully acknowledged. augner@uni-wuppertal.de}}
\title{Well-posedness and Stability of Linear Port-Hamiltonian Systems with Nonlinear Boundary Feedback}

\newtheorem{theorem}{Theorem}[section]
\newtheorem{lemma}[theorem]{Lemma}
\newtheorem{proposition}[theorem]{Proposition}
\newtheorem{corollary}[theorem]{Corollary}
\newtheorem{definition}[theorem]{Definition}

\newtheorem{remark}[theorem]{Remark}
\newtheorem{example}[theorem]{Example}

\newtheorem{assumption}[theorem]{Assumption}

\newcommand{\C}{\mathbb{C}}
\newcommand{\R}{\mathbb{R}}
\newcommand{\diag}{\ \mathrm{diag}}

\renewcommand{\H}{\mathcal{H}}

\newcommand{\N}{\mathbb{N}}
\newcommand{\z}{\zeta}

\renewcommand{\sp}[2]{\left( {#1} \mid {#2} \right)}
\newcommand{\speq}[2]{\left( \left( {#1} \mid {#2} \right) \right)}

\newcommand{\dom}{D}

\newcommand{\K}{\mathbb{K}}
\newcommand{\norm}[1]{\left\Vert{#1}\right\Vert}

\newcommand{\abs}[1]{\left|{#1}\right|}
\newcommand{\A}{\mathcal{A}}
\newcommand{\B}{\mathcal{B}}

\newcommand{\ran}{\operatorname{ran}\,}

\renewcommand{\Re}{\operatorname{Re}\,}

\newcommand{\Lip}{\operatorname{Lip}}

\begin{document}
\maketitle

\pagestyle{myheadings}
\thispagestyle{plain}
\markright{PH SYSTEMS WITH NONLINEAR BOUNDARY FEEDBACK}

 \begin{abstract}
  Boundary feedback stabilisation of linear port-Hamiltonian systems on an interval is considered.
  Generation and stability results already known for linear feedback are extended to nonlinear dissipative feedback, both to static feedback control and dynamic control via an (exponentially stabilising) nonlinear controller.
  A design method for nonlinear controllers of linear port-Hamiltonian systems is introduced.
  As a special case the Euler-Bernoulli beam is considered.
 \end{abstract}



\section{Introduction}
 Within the last years there has been done a lot of research on the stability and stabilisation of wave and beam equations.
 Sufficient conditions which are easy to check for asymptotic or even exponential stabilisation of systems of the abstract \emph{port-Hamiltonian} form
  \begin{equation}
   \frac{\partial x}{\partial t}(t,\z)
    = \sum_{k=0}^N P_k \frac{\partial^k (\H x)}{\partial \z^k}(t,\z),
    \quad
    t \geq 0, \ \z \in (0,1)
    \label{eq:evolution_eqn}
  \end{equation}
 where $x(t,\z) \in \K^d \ (\K = \R$ or $\C)$ via suitable dissipative linear boundary conditions have been given.
 Here, for the case $N=1$ (i.e. in particular wave equation, transport equation and Timoshenko beam equation) we should mention \cite{VillegasEtAl_2009}, the PhD thesis \cite{Villegas_2007} and the monograph \cite{JacobZwart_2012}. Also note the more recent article \cite{Engel_2013} for the case $N = 1$ and $\H = I$.
 More recently, investigations have been done in three generalising directions.
 First, stabilisation using a dynamic controller (\cite{RamirezZwartLeGorrec_2013}, \cite{AugnerJacob_2014}), secondly considering the case $N \geq 2$ (\cite{AugnerJacob_2014}) which includes the Schr\"odinger equation and the Euler-Bernoulli beam equation and last but not least nonlinear feedback (\cite{LeGorrec_2014}, \cite{Trostorff_2014}).
 This article should be seen as a continuation of \cite{AugnerJacob_2014} in which we extend the results given therein for the linear feedback case to the situation of nonlinear feedback.
 Since in \cite{AugnerJacob_2014} the main tool for investigating stability were the Arendt-Batty-Lyubich-V\~u Theorem and the Gearhart-Greiner-Pr\"uss Theorem which both only hold for the case of linear evolution equations we have to develop alternative tools to attack the nonlinear feedback problem.
 Note that the infinite-dimensional system itself remains linear, so we do not touch the topic of nonlinear port-Hamiltonian systems (in the strict sense).
 Still, nonlinear feedback urges us to consider nonlinear contraction semigroups (see, e.g. \cite{Miyadera_1992} and \cite{Showalter_1997}) instead of linear semigroups, following ideas similar to those in \cite{Trostorff_2014} for the generation results and then exploiting ideas already used in \cite{ChenEtAl_1987} (in a linear scenario) for stability properties.
 We should also mention that the approach of \cite{Villegas_2007} (there $N = 1$ and linear feedback) may be used to obtain stability results for both the static and dynamic scenario (also see \cite{LeGorrec_2014}).
 However, most likely this method is restricted to the case $N = 1$.
 Also we stress that the interest in nonlinear feedback is motivated by applications where sometimes (due to physical or technical restrictions) it is not possible to implement perfectly linear controllers.
 Somehow the results of this article therefore show that to some extend nonlinear perturbations (from the perfectly linear case) do not harm the stabilisation properties.
 Also note that in some cases the (usually finite dimensional) control systems considered here actually consist of both a finite dimensional controller and a finite dimensional control target which are connected mechanically via a beam modelled by a infinite-dimensional port-Hamiltonian system, e.g. wave equation, Timoshenko beam or Euler-Bernoulli beam, see e.g. \cite{LeGorrec_2014}.

 This article is structured as follows.
 In Section \ref{sec:background} we recall some background on nonlinear contraction semigroups and $m$-dissipative operators on Hilbert spaces, as well as the Komura-Kato Theorem which as generation theorem for nonlinear contraction semigroups may be seen as the nonlinear analogon to the Lumer-Phillips Theorem.
 Then in Section \ref{sec:impedance_passive_phs} we recall and stress some properties of port-Hamiltonian systems in impedance passive boundary control and observation formulation.
 These observations together with results of Section \ref{sec:background} then give the generation result Theorem \ref{thm:generation-nl-static} (cf. Theorem 5.4 in \cite{Trostorff_2014}) which is restricted to static boundary conditions.
 In fact, we prove that the port-Hamiltonian operator $A$ associated to the evolution equation (\ref{eq:evolution_eqn}) with suitable dissipative boundary conditions generates a (nonlinear) contraction semigroup $(S(t))_{t\geq0}$ on $\overline{\dom(A)} (\ = X := L^2(0,1;\K^d)$ under suitable assumptions).
 From there we first consider the case $N = 1$ and static feedback in Section \ref{sec:N=1-static} which serves as a introductory model case for the more general results later on.
 We are particularly interested in exponential stability of the nonlinear semigroup $(S(t))_{t\geq0}$ (in other words, global exponential stability of the equilibrium $x = 0$), i.e. we ask whether there are constants $M \geq 1$ and $\omega < 0$ such that for all $x \in \overline{\dom(A)} (= X$ in this section) one has
  \begin{equation}
   \norm{S(t) x}_X
    \leq M e^{\omega t} \norm{x}_X,
    \quad t \geq 0.
    \nonumber
  \end{equation}
 We show that under similar conditions as in the linear feedback case, exponential stability can be ensured.
%  whereas in Subsection \ref{subsec:N=1_less_than_linear} we have a look on the situation where the dissipation condition at the boundary is weakened, therefore only obtaining asymptotic stability, i.e. $S(t) x \rightarrow 0$ as $t \rightarrow \infty$ without a uniform decay.
 As example we consider boundary stabilisation of the wave equation.
 In Section \ref{sec:generation-dynamic} we leave the static feedback setting and consider a dissipative interconnection with a nonlinear dynamic controller.
 The main result of that section is the generalisation of Theorem \ref{thm:generation-nl-static} to its dynamic counterpart Theorem \ref{thm:generation_nonlinear_controller}.
 Then the subsequent sections are devoted to a) transferring the results of Section \ref{sec:N=1-static} to the dynamic scenario and b) generalising these results (both the static and the dynamic feedback cases) to port-Hamiltonian systems of order $N = 2$ where we also have a look on the Euler-Bernoulli Beam as a special case (where due to structural assumptions the dissipation assumptions for stabilisation are less restrictive).

\section{Some Background on Contraction Semigroups}
 \label{sec:background}
 Before actually starting with the investigation of port-Hamiltonian systems we first recall some well-known concepts and results on semigroup theory.
 Since we only consider dissipative systems here, we restrict ourselves to the contractive case.
 For the general theory, see \cite{Miyadera_1992}, and for the linear case see, e.g. the monograph \cite{EngelNagel_2000}.
 Throughout we use the following definition.
 \begin{definition}[Semigroup]
  Let $X$ be a Banach space and $X_0 \subset X$ a closed subset.
  A family $(S(t))_{t\geq0}$ of mappings $S(t): X_0 \rightarrow X_0 \ (t \geq 0)$ is called \emph{semigroup} if it satisfies the properties
   \begin{enumerate}
    \item
     $S(0) = I_{X_0}$, the identity map on $X_0$, and
    \item
     $S(t+s) = S(t) S(s)$ for all $s, t \geq 0$.
   \end{enumerate}
  We speak of a \emph{strongly continuous} (abbr.: s.c.) (nonlinear) \emph{semigroup} (or, dynamical system) if $S(t) \in C(X_0;X_0) \ (t \geq 0)$ and for all $x \in X_0$ the map $t \mapsto S(t)x$ is continuous on $\R_+ := [0, \infty)$.
  If additionally $X_0 = X$ and all maps $S(t) \ (t \geq 0)$ are linear, i.e. $S(t) \in \B(X)$, then we speak of a \emph{strongly continuous semigroup} (of linear operators), or \emph{$C_0$-semigroup}.
  A (linear or nonlinear) semigroup $(S(t))_{t\geq0}$ is called \emph{contractive}, if all maps $S(t) \ (t \geq 0)$ are contractions, i.e.
   \begin{equation}
    \norm{S(t)x - S(t)x'}_X
      \leq \norm{x - x'}_X,
      \quad
      x, x' \in X, \ t \geq 0.
      \nonumber
   \end{equation}
 \end{definition}

%  \begin{remark}
%   In the following we denote nonlinear strongly continuous semigroups by $(S(t))_{t\geq0}$, whereas the notation $(T(t))_{t\geq0}$ is reserved for $C_0$-semigroups of linear operators.
%  \end{remark}
 
 Below we restrict ourselves to the Hilbert space case.
 We are going to state the nonlinear version of the Lumer-Phillips Theorem, i.e. the Komura-Kato Theorem, and therefore recall the concepts of dissipative (resp. monotone) operators.
 For details see, e.g. Chapter IV in \cite{Showalter_1997}.
 
 \begin{definition}
  Let $A: X \rightarrow \mathcal{P}(X) := \{ B \subseteq X\}$ be a map.
  We write $\dom(A)$ for its \emph{domain}
   \begin{equation}
    \dom(A)
     := \{x \in X: A(x) \not= \emptyset\}
     \nonumber
   \end{equation}
  and also write $A: \dom(A) \subseteq X \rightrightarrows X$.
  If $A(x) =\{y_x\}$ for all $x \in \dom(A)$ we call $A$ an \emph{operator} and write
   \begin{equation}
   A x
    := y_x,
    \quad
    \text{for the unique} \ y_x \in A(x).
    \nonumber
   \end{equation}
  Otherwise we say that $A$ is \emph{multi-valued}.
 \end{definition}

 We use the notation $A + B$ for the sum of two maps $A: \dom(A) \subset X \rightrightarrows X$ and $B: \dom(B) \subset X \rightrightarrows X$ as follows
  \begin{equation}
   (A + B)(x)
    = \{y_1 + y_2 \in X: y_1 \in A(x), \ y_2 \in B(x) \}.
    \nonumber
  \end{equation}
 Note that $\dom(A + B) = \dom(A) \cap \dom(B)$ and for the particular case where $B$ is an operator $(A + B)(x) = \{y_1 + Bx: y_1 \in A(x)\}$.

 \begin{definition}
  Let $X$ be a Hilbert space and $A: \dom(A) \subseteq X \rightrightarrows X$.
  We say that $A$ is \emph{dissipative} (and $-A$ monotone (or, accretive)), if for all $x, x' \in \dom(A)$ and $y \in A(x), \ y' \in A(x')$ one has
   \begin{equation}
    \Re \sp{y - y'}{x - x'}_X
     \leq 0.
     \nonumber
   \end{equation}
  If additionally
   \begin{equation}
    \{y \in X: \exists x \in \dom(A): y \in (I - A)(x) \}
     =: \ran (I - A)
     = X
     \nonumber
   \end{equation}
  then $A$ (resp. $-A$) is called \emph{m-dissipative} (resp. m-monotone (or, m-accretive)).
  If $A: \dom(A) \subset X \rightrightarrows X$ is dissipative and has no proper dissipative extension, i.e. if $B: \dom(B) \subset X \rightrightarrows X$ with $\dom(A) \subset \dom(B)$ and $A(x) \subset B(x)$ for all $x \in \dom(A)$, then $A = B$, then we call $A$ maximal dissipative.
 \end{definition}

 \begin{remark}
  Let $A: \dom(A) \rightrightarrows X$ be an m-dissipative map on a Hilbert space $X$, then for all $x \in \dom(A)$ the set $A(x)$ is convex and thus there is a unique $z \in A(x)$ with minimal norm.
  This defines the \emph{minimal section} $A^0$ of $A$:
   \begin{equation}
    A^0 x
     := z,
     \quad
     \norm{z}
     = \inf_{y \in A(x)} \norm{y},
     \quad
     \dom(A^0)
     = \dom(A).
     \nonumber
   \end{equation}
%  Also we may define
%   \begin{equation}
%    \abs{A x}
%     := \norm{A^0 x}
%     = \inf_{y \in \dom(A(x))} \norm{y}.
%     \nonumber
%   \end{equation}
  Moreover for all $x \in X$ and $\lambda \in \K_0^+$ the element $y \in \dom(A)$ such that
   \begin{equation}
    x
     \in (\lambda I - A)(y)
     \nonumber
   \end{equation}
  is uniquely determined, thus we may write $y = (\lambda I - A)^{-1}x$.
  In particular every m-dissipative operator is maximal dissipative.
 \end{remark} 
%  \textbf{Proof.}
%  Let $x \in \dom(A)$ be arbitrary and $y, z \in A(x)$, $\lambda \in [0,1]$.
%  Then, since $A$ is m-dissipative there is $\tilde x \in \dom(A)$ such that
%  \begin{equation}
%    \lambda y + (1-\lambda) z - x
%     \in (A-I) (\tilde x),
%      \nonumber
%  \end{equation}
%  i.e. $\lambda y + (1-\lambda) z + \tilde x - x \in A(\tilde x)$.
%  From the dissipativity of $A$ we thus obtain
%   \begin{align}
%    0
%     &\leq \norm{x -\tilde x}^2
%     = \sp{\lambda y + (1-\lambda) z - (\lambda y + (1-\lambda) z - x + \tilde x)}{x - \tilde x}
%     \nonumber \\
%     &= \lambda \sp{y - (\lambda y + (1-\lambda) z - x + \tilde x)}{x - \tilde x}
%     \nonumber \\
%     &\quad + (1 - \lambda) \sp{z - (\lambda y + (1-\lambda) z - x + \tilde x)}{x - \tilde x}
%     \leq 0,
%     \nonumber
%   \end{align}
%  hence $\tilde x = x$ and $\lambda y + (1-\lambda) z \in A(x)$.
% 
%  For the second statement let $\lambda > 0$ and $x \in (A - \lambda I)(y) \cap (A - \lambda I)(z)$,
%  then $x + \lambda y \in A(y)$ and $x + \lambda z \in A(z)$, so
%   \begin{align}
%    0
%     &\leq \sp{y - z}{y - z}
%     = \frac{1}{\lambda} \sp{(x + \lambda y) - (x + \lambda z)}{y - z}
%     \leq 0
%     \nonumber
%   \end{align}
%  and it follows $y = z$.
%  For the last statement note that any dissipative extension $\tilde A$ of an m-dissipative operator is again m-dissipative, thus $\lambda I - \tilde A$ is injective as we just saw, but as extension of the surjective map $\lambda I - A$ both maps then have to be equal, i.e. $\tilde A = A$ is no proper extension.
%  \qed
 
 \begin{lemma}
  \label{lem:perturbation_m-monotone}
  If $A: \dom(A) \rightrightarrows X$ is m-dissipative and $B: X \rightarrow X$ is dissipative and Lipschitz continuous, then also $A + B: \dom(A) \rightrightarrows X$ is m-dissipative.
 \end{lemma}
  \textbf{Proof.}
   See Lemma IV.2.1 in \cite{Showalter_1997}.
   \qed
%  \textbf{Proof.}
%  We follow the line of proof for Lemma IV.2.1 in \cite{Showalter_1997}.
%  First, we note that the sum of two dissipative operators is again dissipative (take the intersection of their domains as the domain of the sum) and also multiplicating an (m-)dissipative operator by $\alpha > 0$ gives another (m-)dissipative operator.
%  Writing $A + B = \frac{1}{\alpha} (\alpha A + \alpha B)$ where $\alpha > 0$ we may therefore assume that $B$ is a strict contraction (i.e. its Lipschitz constant is strictly less than 1).
%  For given $f \in X$ we shall find $x \in \dom(A+B) = \dom(A)$ such that $f \in x - A(x) - Bx$, i.e.
%   \begin{equation}
%    x
%     = \Phi(x)
%     := (I - A)^{-1} (f - Bx).
%     \nonumber
%   \end{equation}
%  Here $\Phi$ is a strict contraction and thus from the strict contraction principle this equation has a unique solution $x \in \dom(A)$.
%  \qed

 For clarity, let us also mention Minty's Theorem.

 \begin{theorem}[Minty]
  On a Hilbert space $X$ the $m$-dissipative operators are exactly the maximal dissipative operators.
 \end{theorem}
 \textbf{Proof.}
 Combine Lemma 2.2.12(iii) and Corollary 3.2.27 in \cite{Miyadera_1992}.
 \qed

 As in the linear ($C_0$-semigroup) case, $m$-dissipative operators are closely related to the generators of contraction semigroups.

 \begin{definition}
  Let $(S(t))_{t\geq0}$ be a nonlinear strongly continuous contraction semigroup on $X$.
  Set
   \begin{equation}
    \hat D
     := \left\{ x \in X: S(\cdot) x \in \Lip(\R_+;X) \right\}.
     \nonumber
   \end{equation}
  We define the (infinitesimal) \emph{generator} of the s.c. contraction semigroup $(S(t))_{t\geq0}$ as
   \begin{equation}
    A_0 (x)
     := \lim_{t \searrow 0} \frac{S(t)x - x}{t},
     \quad
    \dom(A_0)
     := \{x \in X: \lim_{t \searrow 0} \frac{S(t)x - x}{t} \in X \ \text{exists} \}
     \nonumber
   \end{equation}
  and the (g)-operator $A: \dom(A) \subset X \rightrightarrows X$ as the maximal dissipative extension of $A_0$ with $\dom(A) \subset \hat D$.
 \end{definition}

 \begin{remark}
  By Zorn's Lemma every dissipative operator has a maximal dissipative extension (see Lemma 2.2.12(ii) in \cite{Miyadera_1992}).
  Hence the (g)-operator always exists.
  Also note that the infinitesimal operator $A_0$ (or the (g)-operator $A$) uniquely determines the s.c. contraction semigroup (see Corollary 3.4.17 in \cite{Miyadera_1992}).
 \end{remark}
 
%  \begin{theorem}[Lumer-Philipps]
%   Let $A$ be a linear operator on a Hilbert space $X$.
%   Then $A$ generates a contractive $C_0$-semigroup $(T(t))_{t\geq0}$ on $X$ if and only if $A$ is m-dissipative.
%  \end{theorem}
%  \textbf{Proof.}
%  See e.g. Theorem II.3.5 in \cite{EngelNagel_2000}.
%  \qed 
 
 \begin{theorem}[Komura-Kato]
  Let $A: \dom(A) \subseteq X \rightrightarrows X$ be a (possibly multi-valued) map on a Hilbert space $X$.
  If $A$ is m-dissipative, then it generates a nonlinear strongly continuous contraction semigroup $(S(t))_{t\geq0}$ on $\underline X := \overline{\dom(A)}^X$.
  More precisely, for each $x_0 \in \dom(A)$ there is a unique absolutely continuous solution $x \in W^1_\infty(\R_+;X)$ of the abstract nonlinear Cauchy problem
   \begin{align}
    \frac{d}{dt} x(t)
     &\in A(x(t)),
     \quad t \geq 0
     \nonumber \\
    x(0)
     &= x_0.
    \label{eqn:Cauchy_nl} 
   \end{align}
  Also $\norm{\frac{d}{dt} x}_{L_\infty(\R_+;X)} \leq \norm{A^0 x_0}_X$, the function $\norm{A^0 x}_X$ is decreasing and for every $t \geq 0$ and the right-derivative $\frac{d^+}{dt}$ one has
   \begin{equation}
    \frac{d^+}{dt} x(t)
     = A^0 x(t),
     \quad
     t \geq 0.
     \nonumber
   \end{equation}
 \end{theorem}
 \textbf{Proof.}
 See Proposition IV.3.1 in \cite{Showalter_1997}.
 \qed

 \begin{remark}
  If $A$ is m-dissipative and $0 \in A(0)$, then $S(t)(0) = 0$ for all $t \geq 0$.
  Consequently in this case
   \begin{equation}
    \norm{S(t) x}_X
     \leq \norm{x}_X,
     \quad
     t \geq 0.
     \nonumber
   \end{equation}
 \end{remark}

\section{Impedance Passive Port-Hamiltonian Systems}
 \label{sec:impedance_passive_phs}
 In this section we lay the foundations for the generation theorems later on as we introduce port-Hamiltonian systems and boundary control and observation systems.
 We assume impedance passivity (as boundary control and observation system) and observe that the transfer function exists on $\K_0^+$ and its symmetric parts are coercive as linear operators on $\K^{Nd}$.

 \begin{definition}[Port-Hamiltonian System]
  \label{def:PHS}
  Let $N \in \N$ and $P_k \in \K^{d \times d} \ (k = 0, 1, \ldots, N)$ with $P_k^* = (-1)^{k+1} P_k \ (k = 1, \ldots, N)$ and the symmetric part $\Re P_0 := \frac{P_0 + P_0^*}{2} \leq 0$ of the matrix $P_0$ be negative semi-definite as well as $\H \in L_{\infty}(0,1)^{d \times d}$ and consider the Hilbert space $X = L_2(0,1)^d$.
  Further let $W_B, W_C \in \K^{Nd \times 2Nd}$ be matrices such that $\left[ \begin{smallmatrix} W_B \\ W_C \end{smallmatrix} \right] \in \K^{2Nd \times 2Nd}$ is invertible.
  \begin{enumerate}
   \item
    If $P_N$ is invertible and $\H$ is coercive as multiplication operator on $X$, i.e. there is $m_0 > 0$ such that
     \begin{equation}
      z^* \H(\z) z
       \geq m_0 \abs{z}^2,
       \quad
        z \in \K^d, \ \text{a.e.} \ \z \in (0,1)
        \nonumber
     \end{equation}
    then the operator
     \begin{align}
      \mathfrak{A} x
       &= \sum_{k=0}^N P_k \frac{\partial^k}{\partial \z^k} (\H x)
       \nonumber \\
      \dom(\mathfrak{A})
       &= \left\{ x \in L_2(0,1)^d: \ \H x \in H^N(0,1)^d \right\}
       \nonumber
     \end{align}
    is called (maximal) \emph{port-Hamiltonian operator}.
   \item
    For a port-Hamiltonian operator $\mathfrak{A}$ the boundary port-variables $f_{\partial, \H x} \in \K^{Nd}$ (\emph{boundary flow}) and $e_{\partial, \H x} \in \K^{Nd}$ (\emph{boundary effort}) are defined as
     \begin{equation}
      \left( \begin{array}{c} f_{\partial, \H x} \\ e_{\partial, \H x} \end{array} \right)
       = R_{ext} \left( \begin{smallmatrix} (\H x)(1) \\ \vdots \\ (\H x)^{(N-1)}(1) \\ (\H x)(0) \\ \vdots \\ (\H x)^{(N-1)}(0) \end{smallmatrix} \right)
       \nonumber
     \end{equation}
     where $R_{ext} = \frac{1}{\sqrt{2}} \left[ \begin{smallmatrix} Q & -Q \\ I & I \end{smallmatrix} \right] \in \K^{2Nd \times 2Nd}$ for the (invertible) matrix
     \begin{equation}
      Q = \left[ \begin{smallmatrix} P_1 & P_2 & \cdots & \cdots & P_N \\ - P_2 & - P_3 & \cdots & -P_N & 0 \\ \vdots &&& \vdots \\ (-1)^{N-1} P_N & 0 & \cdots & 0 & 0 \end{smallmatrix} \right].
      \nonumber
     \end{equation}
   \item
    For a port-Hamiltonian operator $\mathfrak{A}$ we define the input map $\mathfrak{B}: \dom(\mathfrak{A}) \subset X \rightarrow \K^{Nd}$ and the output map $\mathfrak{C}: \dom(\mathfrak{A}) \subset X \rightarrow \K^{Nd}$ via
     \begin{align}
      \mathfrak{B} x
       &= W_B \left( \begin{smallmatrix} f_{\partial,\H x} \\ e_{\partial, \H x} \end{smallmatrix} \right)
       \nonumber \\
      \mathfrak{C} x
       &= W_C \left( \begin{smallmatrix} e_{\partial,\H x} \\ f_{\partial, \H x} \end{smallmatrix} \right)
       \nonumber
     \end{align}
    and call $\mathfrak{S} = (\mathfrak{A}, \mathfrak{B}, \mathfrak{C})$ a port-Hamiltonian system to which we associate the following abstract boundary control and observation problem
     \begin{align}
      \frac{d}{d t} x(t)
       &= \mathfrak{A} x(t),
       \quad
       x(0) = x_0
       \nonumber \\
      u(t)
       &= \mathfrak{B} x(t)
       \nonumber \\
      y(t)
       &= \mathfrak{C} x(t),
       \quad
       t \geq 0.
       \nonumber
     \end{align}
  \end{enumerate}
 \end{definition}

From here on we assume that $\mathfrak{S} = (\mathfrak{A}, \mathfrak{B}, \mathfrak{C})$ is a port-Hamiltonian system in the sense of Definition \ref{def:PHS}.
Further we always assume that the system $\mathfrak{S} = (\mathfrak{A}, \mathfrak{B}, \mathfrak{C})$ is impedance passive.

 \begin{assumption}
  \label{asmp:passivity}
  The port-Hamiltonian system $\mathfrak{S} = (\mathfrak{A}, \mathfrak{B}, \mathfrak{C})$ is \emph{impedance passive}, i.e.
   \begin{equation}
    \Re \sp{\mathfrak{A} x}{x}_\H
     \leq \Re \sp{\mathfrak{B} x}{\mathfrak{C} x}_{\K^{Nd}},
     \quad
     x \in \dom(\mathfrak{A})
     \nonumber
   \end{equation}
  where we take $X = L_2(0,1)^d$ to be equipped with the inner product $\sp{\cdot}{\cdot}_\H := \sp{\cdot}{\H \cdot}_{L_2}$.
 \end{assumption}

 \begin{remark}
  Assumption \ref{asmp:passivity} already implies that the operators $\mathfrak{A}|_{\ker \mathfrak{B}}$ and $\mathfrak{A}|_{\ker \mathfrak{C}}$ are dissipative on $X$.
  Hence by Theorem 7.2.4 in \cite{JacobZwart_2012} these operators generate s.c. contraction semigroups (of linear operators), respectively.
  In fact, below we prove a nonlinear version of this result, see Theorem \ref{thm:generation-nl-static}.
 \end{remark}

One may even say more, namely $\mathfrak{S} = (\mathfrak{A}, \mathfrak{B}, \mathfrak{C})$ also is a boundary control and observation System.

 \begin{proposition}
  The impedance passive port-Hamiltonian system $\mathfrak{S} = (\mathfrak{A}, \mathfrak{B}, \mathfrak{C})$ is a \emph{boundary control and observation system}, i.e.
   \begin{enumerate}
     \item
      $\mathfrak{A}: \dom(\mathfrak{A}) \subset X \rightarrow X, \ \mathfrak{B}: \dom(\mathfrak{A}) \subset X \rightarrow U, \ \mathfrak{C}: \dom(\mathfrak{A}) \subset X \rightarrow Y$ are linear operators where $X, U$ and $Y$ are Hilbert spaces. (Here: $U = Y = \K^{Nd}$.)
     \item
      $\mathfrak{A}|_{\ker \mathfrak{B}}$ generates a $C_0$-semigroup on $X$.
      (Here: The semigroup is even contractive.)
     \item
      There is a right-inverse $B \in \B(U,X)$ of $\mathfrak{B}$ such that
       \begin{equation}
        \ran B \subseteq \dom(\mathfrak{A}), \
         \mathfrak{A} B \in \B(U,X), \
         \mathfrak{B} B = I.
         \nonumber
       \end{equation}
      \item
       $\mathfrak{C}$ is bounded from $\ker \mathfrak{B}$ to $Y$ where $\ker \mathfrak{B}$ is equipped with the graph norm of $\mathfrak{A}|_{\ker \mathfrak{B}}$.
       (Here: $\mathfrak{C}$ is also bounded from $\dom(\mathfrak{A})$ to $Y$ where $\dom(\mathfrak{A})$ is equipped with the graph norm of the (here: closed) operator $\mathfrak{A}$.)
   \end{enumerate}
 \end{proposition}
 \textbf{Proof.}
 See Theorem 4.4 in \cite{LeGorrecZwartMaschke_2005}.
 \qed

 \begin{remark}
  If $W_B = \left[ \begin{smallmatrix} W_{B,1} \\  W_{B,2} \end{smallmatrix} \right]$ for $W_{B,1} \in \K^{k \times 2Nd}$ and $W_{B,2} \in \K^{(Nd-k) \times 2Nd}$ for $k \in \{0, 1, \ldots, Nd\}$ (and an analogue decomposition of $W_C, \ \mathfrak{B}$ and $\mathfrak{C}$, respectively) and one considers
   \begin{equation}
    \mathfrak{A}_1
     := \mathfrak{A}|_{\dom(\mathfrak{A}_1)},
     \quad
     \dom(\mathfrak{A}_1)
      := \{ x \in \dom(\mathfrak{A}): \ \mathfrak{B}_2 x = 0 \},
      \nonumber
   \end{equation}
  then (quite naturally) also $(\mathfrak{A}_1, \mathfrak{B}_1, \mathfrak{C}_1)$ is a boundary control and observation system (where $U_1 = Y_1 = \K^k$), since this only means that we fix some components of the input $u$ to be zero.
  Hence without loss of generality we always assume $k = Nd$.
  (Here we use the notations $\K^0 := \{0\}$ and $\K^{m \times n} := \B(\K^n,\K^m)$ also in the cases where $n$ or $m$ equals zero.)
 \end{remark}
 
 Let us also introduce the concept of a transfer function which is closely related to the Laplace transform of the semigroup $(T(t))_{t\geq0}$ generated by $A$.
 (For more background on transfer functions for port-Hamiltonian systems we refer to Chapter 12 in \cite{JacobZwart_2012}.)
 
 \begin{definition}[Transfer function]
  Consider the abstract boundary control and observation problem
   \begin{align}
    \frac{d}{dt} x(t)
     &= \mathfrak{A} x(t)
     \nonumber \\
    u(t)
     &= \mathfrak{B} x(t)
     \nonumber \\
    y(t)
     &= \mathfrak{C} x(t),
     \quad
     t \geq 0
     \nonumber
   \end{align}
    where $\mathfrak{S} = (\mathfrak{A}, \mathfrak{B}, \mathfrak{C})$ is a boundary control and observation system
    and let $\lambda \in \K$.
    We write $\lambda \in \dom(G)$ if there is $G(\lambda) \in \B(U,Y)$ such that for all $u \in U$ there is a unique solution of
     \begin{align}
      \lambda x
       &= \mathfrak{A} x
       \nonumber \\
      u
       &= \mathfrak{B} x
       \nonumber \\
      y
       &= \mathfrak{C} x
       \label{eqn:transfer_fct}
     \end{align}
    where $x \in \dom(\mathfrak{A})$ and $y \in Y$ is given by $y = G(\lambda) u$.
 \end{definition}

For impedance passive port-Hamiltonian systems one has $\K_0^+ := \{\lambda \in \K: \Re \lambda > 0 \} \subseteq \dom(G)$.

 \begin{lemma}
  \label{lem:transfer_fct_phs}
  Let $\mathfrak{S} = (\mathfrak{A}, \mathfrak{B}, \mathfrak{C})$ be an impedance passive port-Hamiltonian system with $\H = I$.
  Then $\K_0^+ \subseteq \dom(G)$ and $\Re G(\lambda) > 0$ for all $\lambda \in \K_0^+$, i.e. for all $\lambda \in \K_0^+$ there is $m_\lambda > 0$ such that
   \begin{equation}
    \Re \sp{z}{G(\lambda) z}_U
     > m_{\lambda} \abs{z}^2,
     \quad
     z \in U = \K^{Nd}.
     \nonumber
   \end{equation}
  More precisely, for every $\lambda \in \K_0^+$ there are operators $\Phi(\lambda) \in \B(X), \Psi(\lambda) \in \B(U,X)$ and $F(\lambda) \in \B(X,Y)$ such that for all $f \in X$ and $u \in U$ there is a unique solution of the problem
   \begin{align}
    (\lambda - \mathfrak{A}) x
     &= f
     \nonumber \\
    u
     &= \mathfrak{B} x
     \nonumber \\
    y
     &= \mathfrak{C} x
     \nonumber
   \end{align}
  which is given by
   \begin{align}
    x
     &= \Phi(\lambda) f + \Psi(\lambda) u
     \nonumber \\
    y
     &= F(\lambda) f + G(\lambda) u.
     \nonumber
   \end{align}
 \end{lemma}

 \begin{remark}
  The restriction $\H = I$ is not necessary.
  In fact, for any impedance passive boundary control and observation system $\mathfrak{S} = (\mathfrak{A},\mathfrak{B}, \mathfrak{C})$ (on Hilbert spaces $X$ and $U = Y$) and $P \in \B(X)$ any coercive operator on $X$, also $\mathfrak{S}_P = (\mathfrak{A}P,\mathfrak{B}P, \mathfrak{C}P)$ is an impedance passive boundary control and observation system (on $X_P = X$ equipped with $\sp{\cdot}{\cdot}_{X_P} := \sp{\cdot}{P \cdot}_X$) and the transfer function exists on $\K_0^+$ for $(\mathfrak{A},\mathfrak{B}, \mathfrak{C})$ if and only if it exists on $\K_0^+$ for $\mathfrak{S}_P = (\mathfrak{A}P,\mathfrak{B}P, \mathfrak{C}P)$ (the situation is similar for $\Phi, \Psi$ and $F$ as in Lemma \ref{lem:transfer_fct_phs}).
 \end{remark}

 \textbf{Proof of Lemma \ref{lem:transfer_fct_phs}.}
 Let $\lambda \in \K_0^+$, $u \in \K^{Nd}$ and $f \in X$ be given.
 First, observe that the equation
  \begin{equation}
   (\lambda - \mathfrak{A}) x
    = f
    \nonumber
  \end{equation}
 has the general solution $x = \xi_1$ for $\xi := (x, x', \ldots, x^{(N-1)})$ and $\xi(\z)= e^{\z B_\lambda} \xi(0) + q_f(\z)$ where
  \begin{equation}
   B_\lambda
    = \left[ \begin{smallmatrix} 0&1&0&\cdots&0 \\ 0&0&1&\ddots&\vdots \\ \vdots&&\ddots&\ddots&0 \\ 0&\cdots&\cdots&0&1 \\ \lambda P_N^{-1}-P_N^{-1}P_0&-P_N^{-1}P_1&\cdots&\cdots&-P_N^{-1}P_{N-1} \end{smallmatrix} \right].
   \nonumber
  \end{equation}
  and $q_f(\z) = \int_0^\z e^{(\z - s) B_\lambda} \left( \begin{smallmatrix} 0 \\ \vdots \\ 0 \\ -f(s) \end{smallmatrix} \right) ds$.
 Writing $E_\lambda = e^{B_\lambda}$ input and output may be expressed as
  \begin{align}
   u
    &= W_B R_{ext} \left[ \begin{smallmatrix} E_\lambda \\ I \end{smallmatrix} \right] \xi(0) + W_B R_{ext} \left[ \begin{smallmatrix} q_f(1) \\ 0 \end{smallmatrix} \right],
    \nonumber \\
   y
    &= W_C R_{ext} \left[ \begin{smallmatrix} E_\lambda \\ I \end{smallmatrix} \right] \xi(0) + W_C R_{ext} \left[ \begin{smallmatrix} q_f(1) \\ 0 \end{smallmatrix} \right].
    \nonumber
  \end{align}
 Since the system $(\mathfrak{A}, \mathfrak{B}, \mathfrak{C})$ is impedance passive both the matrices $W_B R_{ext} \left[ \begin{smallmatrix} E_\lambda \\ I \end{smallmatrix} \right]$ and $W_C R_{ext} \left[ \begin{smallmatrix} E_\lambda \\ I \end{smallmatrix} \right]$ are invertible since otherwise (choosing $f = 0$ and $\xi(0)$ in the kernel of one of these matrices) $\lambda \in \K_0^+ \cap \sigma(\mathfrak{A}|_{\ker \mathfrak{B}})$ or $\lambda \in \K_0^+ \cap \sigma(\mathfrak{A}|_{\ker \mathfrak{C}})$, in contradiction to $\mathfrak{A}|_{\ker \mathfrak{B}}$ and $\mathfrak{A}|_{\ker \mathfrak{C}}$ being dissipative.
 As a result, for any given $u \in \K^{Nd}$ and $f \in X$ there is a unique solution $(x,y) \in \dom(\mathfrak{A}) \times \K^{Nd}$ and clearly the map $(f,u) \mapsto (x,y) =: \left[ \begin{smallmatrix} \Phi(\lambda) & \Psi(\lambda) \\ F(\lambda) & G(\lambda) \end{smallmatrix} \right] \left[ \begin{smallmatrix} f \\ u \end{smallmatrix} \right]$ is linear and bounded.
 By the same reasoning one finds (for $f = 0$ fixed) the inverse map $G(\lambda)^{-1}: y \mapsto u$, so that $G(\lambda)$ is bijective.
 Further we have for all $u \in \K^{Nd} \setminus \{0\}$ and the corresponding solution $(x,y) \in \dom(\mathfrak{A}) \times \K^{Nd}$ of (\ref{eqn:transfer_fct}) that
  \begin{equation}
   \Re \sp{u}{G(\lambda)u}_{\K^{Nd}}
    \geq \Re \sp{\mathfrak{A}x}{x}_{L_2}
    = \Re \sp{\lambda x}{x}_{L_2}
    = \Re \lambda \norm{x}_{L_2}^2
    > 0
    \nonumber
  \end{equation}
so that in fact the symmetric part $\Re G(\lambda) > 0$ is strictly positive definite.
\qed

\section{A Generation Theorem for Static Feedback}
To begin with we prove the generation theorem for port-Hamiltonian systems with nonlinear dissipative boundary conditions.
We use a strategy very similar to the linear case, in fact the main differences are twofold:
On the one hand we use the nonlinear generalisation of the Lumer-Phillips Theorem, namely the Komura-Kato Theorem, so that again the proof of the generation result reduces to showing (besides dissipativity which is an assumption) a range condition.
In the linear case it proved convenient (however, not necessary) to reduce the generation theorem to the special case where $\H = I$.
As an additional hurdle, the relevant Lemma 7.2.3 in \cite{JacobZwart_2012} has to be formulated in a nonlinear version, and regarding the proof one should replace any reasoning with the (linear) adjoint by an argument which is suitable for the nonlinear situation.
In fact, this is

 \begin{lemma}
 \label{lem:H=I_enough}
  Let $X$ be a Hilbert space and $A: \dom(A) \subseteq X \rightrightarrows X$ be a dissipative, possibly nonlinear and multivalued, map.
  Further assume that $P \in \B(X)$ is strictly coercive.
  Then if $A -  I$ is surjective, so is $AP - I$. 
 \end{lemma}
 \begin{remark}
  Note that this a very special and simple case of Theorem 2 in \cite{CalvertGustafson_1972}.
  Since the proof of Lemma \ref{lem:H=I_enough} is quite elementary we give it nevertheless.
 \end{remark}
 \textbf{Proof of Lemma \ref{lem:H=I_enough}.}
 As a first step, assume that $\norm{P - I} < \frac{1}{2}$, then from Neumann's series we conclude that $\norm{P^{-1}} \leq \frac{1}{1-\norm{P-I}} < \frac{1}{1 - \frac{1}{2}} = 2$ so that
  \begin{equation}
   \norm{P-I} \norm{P^{-1}}
    =: \rho \in (0,1).
    \nonumber
  \end{equation}
 We shall show that for any given $f \in X$ there is $x \in \dom(AP)$ such that
  \begin{equation}
   (AP - I) (x)
    \ni f
    \nonumber
  \end{equation}
 which is equivalent to solving the problem
  \begin{equation}
   (AP - P) (x)
    \ni f + (I - P)x,
    \nonumber
  \end{equation}
 or, since $(A - I)^{-1}$ exists,
  \begin{equation}
   x
    = \Phi_f(x)
    := P^{-1} (A - I)^{-1} \left( f + (I-P)x \right).
    \nonumber
  \end{equation}
 We show that $\Phi_f: X \rightarrow X$ is a strict contraction and therefore admits a unique fixed point $x_f =: (AP - I)^{-1} f$.
 In fact, we have
  \begin{align}
   &\norm{\Phi_f(x) - \Phi_f(x')}
    \nonumber \\
    &\leq \norm{P^{-1}} \norm{(A - I)^{-1}(f + (I-P)x) - (A - I)^{-1}(f + (I-P)x')}
    \nonumber \\
    &\leq \norm{P^{-1}} \norm{(f + (I-P)x) - (f + (I-P)x')}
    \nonumber \\
    &\leq \norm{P^{-1}} \norm{I-P} \norm{x - x'}
    = \rho \norm{x - x'}
    \nonumber
  \end{align}
 where we used Corollary 1.3(b) in \cite{Showalter_1997} in the second step.
 Therefore $\Phi_f$ is a strict contraction and the Contraction Principle gives a unique solution $x_f =: (AP - I)^{-1} f$.
 In the second step we remove the restriction on $P$.
 Namely it is easy to see that there are a number $n \in \N$ and a coercive operator $Q = P^{1/n} \in \B(X)$ such that $\norm{I - Q} < \frac{1}{2}$ and $P = Q^n$.
 Note that $\norm{I - Q} = \norm{I - Q}_k$ for all the norms
  \begin{equation}
   \norm{\cdot}_k
    := \norm{Q^{k/2} \cdot},
    \qquad k = 0, 1, \ldots, n-1.
    \nonumber
  \end{equation}
 Writing
  \begin{equation}
   A P - I
    = (AQ^{n-1})Q - I
    \nonumber
  \end{equation}
 the general case follows by induction using the spaces $X_k := \left( X, \norm{\cdot}_k \right)$.
 \qed

It is an easy consequence of the preceding lemma that for any m-dissipative operator $A: \dom(A) \subseteq X \rightrightarrows X$ and a strictly coercive operator $P \in \B(X)$ also the operator $AP$ with domain $\dom(AP) = \{x \in X: Px \in \dom(A) \}$ is m-dissipative on $X$ equipped with the equivalent inner product $\sp{\cdot}{\cdot}_P = \sp{\cdot}{P\cdot}$ (cf. Lemma 5.1 in \cite{Trostorff_2014}).
Of course, in our particular situation $P = \H$ is the Hamiltonian density multiplication operator.

 \begin{theorem}
  \label{thm:generation-nl-static}
  Let $\mathfrak{S} = (\mathfrak{A}, \mathfrak{B}, \mathfrak{C})$ be an impedance passive port-Hamiltonian system.
  Assume that $\phi: \R^{Nd} \rightrightarrows \R^{Nd}$ is a (possibly multi-valued, nonlinear) m-monotone map.
  Then the (single-valued) operator
   \begin{align}
    A
     &= \mathfrak{A}|_{\dom(A)}
     \nonumber \\
     \dom(A)
     &= \left\{ x \in \dom(\mathfrak{A}): \ \mathfrak{B} x \in - \phi(\mathfrak{C} x) \right\}
     \nonumber
   \end{align}
  generates a s.c. contraction semigroup on $X = L_2(0,1)^d$ with the inner product $\sp{\cdot}{\cdot}_{\H} = \sp{\cdot}{\H \cdot}_{L_2}$.
 \begin{remark}
  Note that for the case $N = 1$ a characterisation of $m$-dissipative boundary conditions yielding an $m$-dissipative operator $A$ has been given in Theorem 5.4 of \cite{Trostorff_2014}.
  Also note the more general result Theorem 3.1 therein.
 \end{remark}
 \end{theorem}
 \textbf{Proof of Theorem \ref{thm:generation-nl-static}.}
 From Lemma \ref{lem:H=I_enough} we know that it suffices to consider the case where $\H = I$ equals the identity.
 Also note that there is $x_0 \in \dom(A) \not= \emptyset$ which implies that $x_0 + C_c^\infty(0,1)^d \subseteq \dom(A)$ is a dense subset of $X$.
 Clearly $A$ is dissipative since for $x, \tilde x \in \dom(A)$ we have
  \begin{align}
   \Re \sp{A(x) - A(\tilde x)}{x - \tilde x}_{L_2}
    &= \Re \sp{\mathfrak{A}(x - \tilde x)}{x - \tilde x}_{L_2}
    \nonumber \\
    &\leq \Re \sp{\mathfrak{B}x-\mathfrak{B}\tilde x}{\mathfrak{C}x-\mathfrak{C} \tilde x}_{\K^{Nd}}
    \leq 0
    \nonumber
  \end{align}
 using that $\mathfrak{B} x \in - \phi(\mathfrak{C}x), \ \mathfrak{B}\tilde x \in - \phi(\mathfrak{C}\tilde x)$ and $\phi$ is monotone.
 It remains to show that $\ran (I - A) = X$, i.e. for every $f \in X$ we have to find $x \in \dom(\mathfrak{A})$ such that
  \begin{align}
   (I - \mathfrak{A})x
    &= f
    \nonumber \\
   \mathfrak{B} x
    &\in - \phi(\mathfrak{C} x).
    \nonumber
  \end{align}
 From Lemma \ref{lem:transfer_fct_phs} we know that all solutions of the first of these equations have the form $x = \xi_1$ where
  \begin{equation}
   \xi(\z)
    = e^{\z B_1} \xi(0)
     + \int_0^\z e^{(\z-s)B_1} \left( \begin{smallmatrix} 0 \\ \vdots \\ 0 \\ -f(s) \end{smallmatrix} \right) ds
    \nonumber
  \end{equation}
 and the problem thus reduces to finding $\xi(0) \in \K^{Nd}$ such that
  \begin{equation}
   W_B R_{ext} \left[ \begin{smallmatrix} E \\ I \end{smallmatrix} \right] \xi(0)
    \in - \phi \left( W_C R_{ext} \left[ \begin{smallmatrix} E \\ I \end{smallmatrix} \right] \xi(0) \right)
    \label{thm:generation-nl-static:1}
  \end{equation}
 where $E = E_1 = e^{B_1}$, or, by Lemma \ref{lem:transfer_fct_phs}, finding $y \in \K^{Nd}$ such that
  \begin{equation}
   u
    = G(1)^{-1} y - G(1)^{-1} F(1) f
    \in - \phi(y),
    \nonumber
  \end{equation}
 i.e. $(G(1)^{-1} + \phi)(y) \ni G(1)^{-1} F(1)f$.
 Since $\phi$ is m-monotone and $\Re G(1)^{-1}$ is coercive by Lemma \ref{lem:transfer_fct_phs}, also $\phi + G(1)^{-1} - \varepsilon I$ is m-monotone by Lemma \ref{lem:perturbation_m-monotone} for some small $\varepsilon > 0$.
 We conclude that there is a (unique) $y \in \K^{Nd}$ such that for $u := G(1)^{-1}y$ one has $u \in - \phi(y)$ and since the matrix $W_B R_{ext} \left[ \begin{smallmatrix} E \\ I \end{smallmatrix} \right]$ is invertible it follows that there is a (unique) $\xi(0) \in \K^{Nd}$ such that (\ref{thm:generation-nl-static:1}) holds, i.e. there is a (unique) $x \in \dom(A)$ with $f \in (I - A)(x)$.
 Now the assertion follows from the Komura-Kato Theorem.
 \qed

\section{Stabilisation of First Order Systems}
 \label{sec:N=1-static}
 For the moment we additionally assume that $N = 1$, i.e. $\mathfrak{A} = P_1 (\H \cdot)' + P_0 (\H \cdot)$ on $\dom(\mathfrak{A}) = \{ x \in L_2(0,1)^d: \H x \in H^1(0,1)^d \}$.
 Also we assume that $\H \in W_\infty^1(0,1)^{d \times d}$ is Lipschitz continuous.

We aim to prove the following uniform exponential stability result.
 \begin{theorem}
  \label{thm:exp_stability_nonlinear}
  Let $\mathfrak{S} = (\mathfrak{A}, \mathfrak{B}, \mathfrak{C})$ be an impedance passive port-Hamiltonian system and $\phi: \K^d \rightrightarrows \K^d$ an m-monotone map with $0 \in \phi(0)$.
  For the nonlinear operator
   \begin{equation}
    A
     := \mathfrak{A}|_{\dom(A)},
     \quad
     \dom(A)
     := \{ x \in \dom(\mathfrak{A}): \mathfrak{B} x \in - \phi(\mathfrak{C} x) \}
     \nonumber 
   \end{equation}
  assume that there is $\kappa > 0$ such that
   \begin{equation}
    \Re \sp{A x}{x}_\H
     \leq - \kappa (x^* \H x)(1),
     \quad
     x \in \dom(A).
     \nonumber
   \end{equation}
  Then $A$ generates a s.c. contraction semigroup $(S(t))_{t\geq0}$ with globally exponentially stable equilibrium $0$, i.e. there are $M \geq 1$ and $\omega < 0$ such that
   \begin{equation}
    \norm{S(t) x}_\H
     \leq M e^{\omega t} \norm{x}_\H,
     \quad x \in X, \ t \geq 0.
     \nonumber
   \end{equation}
 \end{theorem}
 \begin{remark}
  If $\phi \in \B(U) = \K^{d \times d}$ is linear this is Theorem III.2 in \cite{VillegasEtAl_2009} which uses a ``sideways energy estimate" (Lemma III.1 in \cite{VillegasEtAl_2009}) in the spirit of \cite{CoxZuazua_1995}.
  Actually a first result like this may already be found as Theorem 3 in \cite{RauchTaylor_1974} where $\H$ is assumed to be smooth.
  An alternative proof of Theorem III.2 in \cite{VillegasEtAl_2009} via Gearhart's Theorem can be found as Proposition 2.12 in \cite{AugnerJacob_2014}, but clearly the latter technique is restricted to the linear situation.
  In fact, we will use a technique which for the linear case yields a third proof of the theorem.
  Namely we use an idea which was used in \cite{ChenEtAl_1987} to prove exponential stability for a chain of linear Euler-Bernoulli beams with (linear) dissipative linkage.
 \end{remark}

 As preparation we state the following auxiliary result.
  \begin{lemma}
  \label{lem:multiplier}
   Let $\alpha > 0$ and $\beta, \gamma \geq 0$ be given.
   Then there is $\eta \in C^{\infty}([0,1];\R)$ with $\eta(0) = 0$ and $\eta' > 0$ such that
    \begin{equation}
     \alpha \eta'(\z) - \beta \eta(\z)
      \geq \gamma,
      \quad
      \z \in [0,1].
      \label{eqn:lem:choice_of_eta}
    \end{equation}
  \end{lemma}
  \textbf{Proof.}
   Scaling $\eta$ by the factor $\frac{1}{\gamma}$ it is enough to consider the case $\gamma = 1$.
   We make the ansatz $\eta(\z) = e^{\lambda \z} - 1$ for $\lambda > 0$ which we are going to specify.
   Then equation (\ref{eqn:lem:choice_of_eta}) is equivalent to
    \begin{align}
     (\alpha \lambda - \beta) e^{\lambda \z}
      \geq \gamma - \beta \ (\z \in [0,1])
      \ &\Leftrightarrow \
      \alpha \lambda  \geq \gamma
      \nonumber \\
      & \Leftrightarrow \
      \lambda \geq \frac{\gamma}{\alpha}.
      \nonumber 
    \end{align}
  \qed
%  \begin{lemma}
%   Let $N \in \N, \ c > 0$ and $M \geq 0$ be given.
%   There is $f \in C^\infty([0,1];\R)$ such that $f(0) = 0, \ f' > 0$ and
%    \begin{equation}
%     f'(\z)
%      > c \norm{f^{(j)}}_{L_\infty} + M,
%      \quad
%      j = 0, \ldots, N, j \not=1, \ \z \in [0,1].
%      \label{eq:star}
%    \end{equation}
%  \end{lemma}
%  \textbf{Proof.}
%  Scaling $f$ by a suitable large factor we may and will assume $M = 0$.
%  Take $f(\z) := e^{\alpha \z} - 1$ for some $\alpha > 0$.
%  Then $f^{(j)}(\z) = \alpha^j e^{\alpha \z} \ (j = 1, \ldots, N)$, so for all $\z, \tilde \z \in [0,1]$
%   \begin{equation}
%    f'(\z)
%     = \alpha^{1-j} e^{\alpha(\z - \tilde \z)} \abs{f^{(j)}(\tilde \z)}
%     \geq \alpha^{1-j} e^{- \alpha} \norm{f^{(j)}}_{L_\infty}
%     \nonumber
%   \end{equation}
%  and there is $\alpha_0 > 0$ such that $\alpha^{1-j} e^{-\alpha} > c \ (j = 2, \ldots, N, \ \alpha \in (0,\alpha_0))$.
%  Secondly
%   \begin{equation}
%    \frac{\abs{f(\tilde \z)}}{f'(\z)}
%     = \frac{e^{\alpha \tilde \z} - 1}{\alpha e^{\alpha \z}}
%     = \frac{e^{\alpha \tilde \z} - 1}{\alpha \tilde \z} \frac{\tilde \z}{e^{\alpha \z}}
%     \rightarrow 0
%     \nonumber
%   \end{equation}
%  as $\alpha \rightarrow 0$, uniformly in $\z, \tilde \z \in (0,1)$, i.e. $f' > c \norm{f}_{L_\infty}$ for all $\alpha \in (0,\alpha_1)$ for some $\alpha_1 > 0$.
%  Finally, if $\alpha \in (0, \min \{\alpha_0, \alpha_1\})$, the assertion follows.
% \qed

 Also the following fact (which can be derived via integration by parts) will prove quite useful in the computations below.
  \begin{lemma}
   \label{lem:real_part}
   Let $Q \in W_\infty^1(0,1)^{d \times d}$ be a function of self-adjoint operators and $u \in H^1(0,1)^d$.
   Then
    \begin{equation}
     \Re \sp{u'}{Qu}_{L_2}
      = - \frac{1}{2} \sp{u}{Q'u}_{L_2} + \frac{1}{2} \left[u(\z)^* Q(\z) u(\z) \right]_0^1.
      \nonumber
    \end{equation}
  \end{lemma}

 \textbf{Proof of Theorem \ref{thm:exp_stability_nonlinear}.}
 Existence of the s.c. contraction semigroup $(S(t))_{t \geq 0}$ follows from Theorem \ref{thm:generation-nl-static}.
 Stability:
 On $X = L_2(0,1)^d$ define the quadratic functional
  \begin{equation}
   q(x)
    := \sp{x}{\eta P_1^{-1} x}_{L_2},
    \quad x \in X
    \nonumber
  \end{equation}
 where $\eta \in C^1([0,1];\R)$ is a function with $\eta(0) = 0$ and $\eta' > 0$ which we choose at a later point.
 Let $x_0 \in \dom(A)$ be arbitrary and denote by $x(t,\cdot) = x(t) = S(t)x_0$ the solution for the initial value $x_0$, so that $x \in W^1_{\infty}(\R_+;X)$, hence $q(x) \in W^1_{\infty}(\R_+;\R)$ and we calculate (using Lemma \ref{lem:real_part}) for a.e. $t \geq 0$
  \begin{align}
   \frac{d}{dt} q(x(t))
    &= 2 \Re \sp{P_1^{-1} \dot x(t)}{\eta x(t)}_{L_2}
    \nonumber \\
    &= 2 \Re \sp{(\H x(t))' + P_1^{-1}P_0(\H x(t))}{\eta x(t)}_{L_2}
     \nonumber \\
    &= - \sp{\H x(t)}{(\eta'\H^{-1} + \eta (\H^{-1})' - 2\eta \Re (\H^{-1} P_1^{-1} P_0))\H x(t)}_{L_2}
     \nonumber \\
     &\quad
     + \left[ \H x(t,\z)^* (\eta \H^{-1})(\z) \H x(t,\z) \right]_0^1.
     \nonumber
  \end{align}
 Thus for $\Phi(t) := t \norm{x(t)}_\H^2 + q(x(t))$ we have
  \begin{align}
   \frac{d}{dt} \Phi(t)
    &= \sp{\H x(t)}{\left( (1-\eta')\H^{-1} + \eta \left[-(\H^{-1})' + 2 \Re \left( \H^{-1} P_1^{-1} P_0 \right) \right] \right) \H x(t)}_{L_2}
    \nonumber \\
    &\quad
     + \eta(1) x(t,1)^* \H(1) x(t,1)
     + 2t \Re \sp{A x(t)}{x(t)}_\H
     \nonumber \\
    &\leq \sp{\H x(t)}{\left( (1-\eta')\H^{-1} + \eta \left[ -(\H^{-1})' + 2 \Re \left( \H^{-1} P_1^{-1} P_0 \right) \right] \right) \H x(t)}_{L_2}
     \nonumber \\
    &\quad
     + (\eta(1) - 2t\kappa) x^*(t,1) \H(1) x(t,1).
    \nonumber
  \end{align}
 So far, we did not specify our choice of $\eta$, so we may choose $\eta$ in such a way that
  \begin{align}
   &(1-\eta'(\z))\H^{-1}(\z)
    + \eta(\z) \left[ (\H^{-1})'(\z) + 2\Re \left( \H^{-1}(\z) P_1^{-1} P_0 \right) \right]
    \nonumber \\
    &\quad \leq \left( M_0 - \eta'(\z) m_0 + \eta(\z) \left[M_1 + 2M_3 \right] \right) I
    \leq 0,
    \quad
    \text{a.e.} \ \z \in (0,1)
    \nonumber
  \end{align}
 where
  \begin{align}
   m_0 I \leq \H^{-1}(\z) \leq M_0 I, \ (\H^{-1})'(\z) \leq M_1 I
    \nonumber \\
   - M_3 I \leq \Re \left(\H^{-1}(\z) P_1^{-1} P_0 \right) \leq M_3 I
   \nonumber
  \end{align}
 for a.e. $\z \in (0,1)$ and we applied Lemma \ref{lem:multiplier}.
 Then for $t \geq t_0 := \frac{\eta(1)}{2\kappa}$ we have
  \begin{equation}
   \frac{d}{dt} \Phi(t)
    \leq 0
    \nonumber
  \end{equation}
 and thus $\Phi$ decreases on $(t_0,\infty)$.
 (Note that the choice of $\eta$ and $t_0$ is independent of the initial value $x_0 \in \dom(A)$.)
 Moreover, since $\abs{q(y)} \leq c \norm{y}^2_\H$ for some $c > 0$ and all $y \in X$ we obtain for $t \geq t_0$ the estimate
  \begin{align}
   t \norm{x(t)}^2_\H
    &\leq \Phi(t) + c \norm{x(t)}^2_\H
    \nonumber \\
    &\leq \Phi(t_0) + c \norm{x(t)}^2_\H
    \nonumber
  \end{align}
 and hence for $t > \max \{t_0, c\}$
  \begin{align}
   \norm{x(t)}^2_\H
    &\leq \frac{\Phi(t_0)}{t - c}
    \leq \frac{t_0 + c}{t - c} \norm{x(t_0)}^2_\H
    \nonumber \\
    &\leq \frac{t_0 + c}{t - c} \norm{x_0}^2_\H
    \nonumber 
  \end{align}
 and from the density of $\dom(A)$ in $X$ we conclude for $t > \max \{t_0, c\}$
  \begin{equation}
   \norm{S(t) x}_X
    \leq \sqrt{\frac{t_0 + c}{t - c}}\norm{x}_X
    \xrightarrow{t \rightarrow + \infty} 0,
    \quad
    x \in X.
    \nonumber
  \end{equation}
 As in the linear case, this property easily implies uniform exponential energy decay.
 \qed

 \begin{remark}
  For $q$ as in the proof of Theorem \ref{thm:exp_stability_nonlinear} the following holds.
  For every solution $x \in W_\infty^1(\R_+;X) \cap L_\infty(\R_+;\dom(\mathfrak{A}))$ of $\dot x = \mathfrak{A} x$ one has $q(x) \in W_\infty^1(\R_+;\R)$ with
   \begin{equation}
    \norm{x(t)}_\H^2 + \frac{d}{dt} q(x(t))
     \leq c \abs{(\H x)(t,1)}^2,
     \quad \text{a.e.} \ t \geq 0.
     \nonumber
   \end{equation}
   We come back to this property in the context of dynamic controllers.
 \end{remark}

 \begin{remark}
  An alternative proof can be established via the ``Sideways Energy Estimate" of Lemma 9.1.2 in \cite{JacobZwart_2012}.
  In fact, the proof of Theorem 9.1.3 in \cite{JacobZwart_2012} almost literally carries over to the nonlinear situation.
 \end{remark}
%  \textbf{``Sideways energy estimate"-based proof of Theorem \ref{thm:exp_stability_nonlinear}.}
%   The proof of Lemma 9.1.2 in \cite{JacobZwart_2012} extends straight-forward to the situation with nonlinear boundary feedback.
%   Thus there are constants $c > 0$ and $\tau > 0$ such that for every $x_0 \in \dom(A)$ and $x := S(\cdot) x_0 \in W_\infty^1(0,\infty;X)$ the estimate
%    \begin{equation}
%      \norm{x(\tau)}_X^2
%       \leq c \int_0^\tau (x^* \H x)(t,1) dt
%       \nonumber
%    \end{equation}
%   holds true.
%   Then
%    \begin{align}
%     \norm{x(\tau)}_\H^2 - \norm{x(0)}_\H^2
%      &= \int_0^\tau \Re \sp{Ax(t)}{x(t)}_\H dt
%      \nonumber \\
%      &\leq - \kappa \int_0^\tau (x^* \H x)(t,1) dt
%      \nonumber \\
%      &\leq - \frac{\kappa}{c} \norm{x(\tau)}_\H^2
%      \nonumber
%    \end{align}
%   and hence $\norm{S(\tau) x_0}_\H \leq \sqrt{\frac{c}{c+k}} \norm{x_0}$ and from time invariance of the problem and the semigroup property it follows
%    \begin{align}
%     \norm{S(t) x_0}_\H
%      &= \norm{S(\tau)^{\lfloor \frac{t}{\tau}\rfloor} S(t - \tau \lfloor \frac{t}{\tau}\rfloor) x_0}_\H
%      \nonumber \\
%      &\leq \left( \sqrt{\frac{c}{c+k}} \right)^{\lfloor \frac{t}{\tau}\rfloor} \norm{x_0}_\H
%      \nonumber \\
%      &\leq \sqrt{\frac{c+k}{c}} e^{- \frac{t}{2\tau} \ln(\frac{c+k}{c})} \norm{x_0}_\H.
%      \nonumber
%    \end{align}
%   \qed

Let us also state an asymptotic stability result which follows from Theorem \ref{thm:exp_stability_nonlinear}, the contraction property and the following interpolation inequality.

 \begin{lemma}
 \label{lem:gn-inequality-1d}
  Let $1 \leq j < N \in \N$ be natural numbers, then there is a constant $C > 0$ such that for every $f \in H^N(0,1)$ one has
   \begin{equation}
    \norm{f^{(j)}}_{C[0,1]}
     \leq C \left( \norm{f}_{L_2(0,1)}^2 + \norm{f^{(N)}}_{L_2(0,1)}^2 \right)^{\frac{2j+1}{2N}} \norm{f}_{L_2(0,1)}^{\frac{2(N-j)-1}{2N}}.
     \nonumber
   \end{equation}
 \end{lemma}
  \textbf{Proof.}
   Combine Lemmas 4.10 and 4.12 in \cite{Adams_1975} to get a one-dimensional version of Theorem 4.14 therein.
   \qed
 \begin{remark}
  Actually Lemma 5.5 together with Theorem IV.1.1 in \cite{Kato_1966} implies that $\mathfrak{A}$ is a closed operator.
 \end{remark}

 \begin{corollary}
  Under the assumptions of Theorem \ref{thm:exp_stability_nonlinear}, but with the less restrictive condition
   \begin{equation}
    \Re \sp{A x}{x}_\H
     \leq - \kappa (x^* \H x)(1),
     \quad
     x \in \dom(A): \
     \abs{\mathfrak{B} x}, \abs{\mathfrak{C} x} \leq \rho
     \nonumber
   \end{equation}
  for some $\rho > 0$, the semigroup $(S(t))_{t\geq0}$ is asymptotically stable, i.e. for all $x \in X$ one has $S(t) x \xrightarrow{t \rightarrow \infty} 0$.
 \end{corollary}
 \textbf{Proof.}
 \emph{First step:} Take any $x_0 \in \dom(A)$ and set $x(t) := S(t)x_0$.
 Then $\norm{x(t)}_{\H}$ and $\norm{\mathfrak{A}x(t)}_{\H}$ are bounded by $\norm{x_0}_{\H}$ and $\norm{\mathfrak{A}x_0}_\H$ for a.e. $t \geq 0$, respectively.
 Hence also $\abs{\mathfrak{B}x(t)}, \abs{\mathfrak{C} x(t)} \leq c \norm{x_0}_{\mathfrak{A}} := c \sqrt{\norm{x_0}_\H^2 + \norm{\mathfrak{A} x_0}_\H^2}$ for a.e. $t \geq 0$, so that $x(t) = \hat S(t) x_0$ where $\hat S(t)$ is the s.c. contraction semigroup corresponding to $\hat A := \mathfrak{A}|_{\dom(\hat A)}$ where $\dom(\hat A) = \{x \in \dom(\mathfrak{A}): \mathfrak{B} x \in - \hat \phi(\mathfrak{C} x) \}$ for some $m$-dissipative $\hat \phi: \dom(\hat \phi) \subseteq \K^{Nd} \rightrightarrows \K^{Nd}$ such that $\hat \phi = \phi$ for $\abs{z} \leq c \norm{x_0}_{\mathfrak{A}}$ and such that
  \begin{equation}
   \Re \sp{\hat A x}{x}_\H
    \leq - \hat \kappa (x^* \H x)(1),
    \quad
    x \in \dom(\hat A).
    \nonumber
  \end{equation}
 Consequently $x(t) = \hat S(t) x_0 \xrightarrow{t \rightarrow \infty} 0$ due to Theorem \ref{thm:exp_stability_nonlinear}. \\
 \emph{Second step:} Let $x \in X$ be arbitrary.
 We have to prove that for all $\varepsilon > 0$ there is $T_\varepsilon > 0$ such that $\norm{S(t)x} \leq \varepsilon$ for all $t \geq T_\varepsilon$.
 For this end take any $x_\varepsilon \in \dom(A)$ such that $\norm{x-x_\varepsilon} \leq \frac{\varepsilon}{2}$.
 Then by the first step for $T_\varepsilon > 0$ sufficiently large one has $\norm{S(t)x_\varepsilon} \leq \frac{\varepsilon}{2}$ \ ($t \geq T_\varepsilon$), so that
  \begin{align}
   \norm{S(t)x}
    &\leq \norm{S(t)x - S(t)x_\varepsilon} + \norm{S(t)x_\varepsilon}
    \nonumber \\
    &\leq \norm{x - x_\varepsilon} + \norm{S(t)x_\varepsilon}
    \leq \varepsilon,
    \quad
    t \geq T_\varepsilon.
    \nonumber
  \end{align}
 \qed

 \begin{example}[Wave Equation]
  Consider the one-dimensional wave equation
   \begin{align}
    \rho \omega_{tt}(t,\z) - (EI \omega_\z)_\z(t,\z)
     &= 0,
     \quad
     \z \in (0,1), \ t \geq 0
     \nonumber
   \end{align}
  where $EI, \rho \in L_\infty(0,1)$ are uniformly positive, i.e. also $EI^{-1}, \rho^{-1} \in L_\infty(0,1)$.
  At the left end we assume conservative or dissipative boundary conditions of the form
   \begin{equation}
    \omega_t(t,0) = 0
     \quad \text{or} \quad
    (EI \omega_\z)(t,0) \in f(\omega_t(t,0)),
     \quad
     t \geq 0
    \nonumber
   \end{equation}
  where $f: \R \rightrightarrows \R$ is maximal monotone and $f(0) \ni 0$, e.g. $f$ could be single-valued, continuous and non decreasing with $f(0) = 0$, in particular the case $f = 0$ (Neumann-boundary condition) is allowed.
  We further assume that on the right end a (monotone) damper is attached to the system, so that the boundary condition is given by
   \begin{equation}
    (EI \omega_\z)(t,1)
     \in - g(\omega_t(t,1))
     \nonumber
   \end{equation}
  where again $g: \R \rightrightarrows \R$ is maximal monotone with $g(0) \ni 0$.
  Of course, the choice $f = g = 0$ would lead to Neumann-boundary conditions on both sides for which the system is known to be energy-preserving, in particular not strongly stable.
  Here as usual the energy is given by
   \begin{equation}
    E(t)
     := \int_0^1 \rho(\z) \abs{\omega_t(t,\z)}^2 + EI(\z) \abs{\omega_\z(t,\z)}^2 d\z.
     \nonumber
   \end{equation}
  In fact the example fits into our port-Hamiltonian setting when we choose $x = (\rho \omega_t, \omega_z)$, $\H = \diag(\rho^{-1}, EI)$ and $P_1 = \left( \begin{smallmatrix} &1\\1& \end{smallmatrix} \right), \ P_0 = 0$ and $\K = \R$.
  If we define
   \begin{align}
    \mathfrak{B} x
     &= \left( \begin{array}{c} - (EI \omega_\z)(0) \\ (EI \omega_\z)(1) \end{array} \right)
     = \left( \begin{array}{c} - (\H x)_2(0) \\ (\H x)_2(1) \end{array} \right)
     \nonumber \\
    \mathfrak{C} x
     &= \left( \begin{array}{c} \omega_t(0) \\ \omega_t(1) \end{array} \right)
     = \left( \begin{array}{c} (\H x)_1(0) \\ (\H x)_1(1) \end{array} \right)
     \nonumber
   \end{align}
  then the system $\mathfrak{S} = (\mathfrak{A}, \mathfrak{B}, \mathfrak{C})$ is impedance passive, since for the maximal operator $\mathfrak{A}$ one has
    \begin{equation}
     \sp{\mathfrak{A} x}{x}_{\H}
      = (\H x)_1(1)^* (\H x)_2(1) - (\H x)_1(0)^* (\H x)_2(0).
      \nonumber
    \end{equation}
  (Note that in the Dirichlet case $\omega_t(0) = 0$ one has to exchange the first components of $\mathfrak{B}$ and $\mathfrak{C}$ and then choose $f = 0$.)
  The corresponding port-Hamiltonian operator $A = \mathfrak{A}|_{\dom(A)}$ (with nonlinear boundary conditions) is dissipative then
   \begin{align}
    \dom(A)
     &= \{x \in L_2(0,1;\R^2): \ \H x \in H^1(0,1;\R^2), \ (\H x)_2(1) \in - g((\H x)_1(1)),
      \nonumber \\
      &\quad
      \begin{cases} (\H x)_1(0) = 0,& \text{(Dirichlet b.c.) \ or} \\ (\H x)_2(0) \in f((\H x)_1(0))& \text{(Neumann b.c.)} \end{cases} \ \}
     \nonumber
   \end{align}
  and we have at least
   \begin{equation}
    \sp{A x}{x}_{\H}
     \leq - (\H x)_1(1) g^0((\H x)_1(1)),
     \quad x \in \dom(A).
     \nonumber
   \end{equation}
  Theorem \ref{thm:generation-nl-static} assures that $A$ generates a nonlinear contraction semigroup on $X = L_2(0,1;\R^2)$ with inner product $\sp{\cdot}{\cdot}_\H$.
  To have stability results we need stronger assumptions on the damper, i.e. on the map $g$.
   First assume (additionally to $g$ being m-monotone) that there is $\kappa > 0$ such that $\kappa^{-1} \abs{x} \leq \abs{z} \leq \kappa \abs{x}$ for all $x \in \R$ and $z \in g(x)$ (i.e. in particular $g(0) = \{0\}$).
   Then we obtain the dissipativity condition
    \begin{equation}
     \sp{Ax}{x}_\H
      \leq - \tilde \kappa \abs{(\H x)(1)}^2,
      \quad x \in \dom(A)
      \nonumber
    \end{equation}
  where $\tilde \kappa := \frac{1}{2} \min \{\kappa, \kappa^{-1}\}$ and so Theorem \ref{thm:exp_stability_nonlinear} ensures uniform exponential stability of the corresponding nonlinear semigroup.
  Secondly we assume that the condition $\kappa^{-1} \abs{x} \leq \abs{z} \leq \kappa \abs{x}$ only holds for $x \in \R$ and $z \in g(x)$ whenever $\abs{x} \leq \rho$ for some fixed $\rho > 0$.
  Then we only obtain asymptotic stability of all solutions $x(t) = S(t) x_0 \ (t \geq 0)$.
  We refer to Example 3.3 in \cite{ConradLeblondMarmorat_1990} for sufficient conditions leading to strong stability of the $n$-dimensional wave equation on a smooth, bounded domain $\Omega \subseteq \R^n$.
 \end{example}
 
 \begin{example}[Dissipation near $0$]
  Let $\phi: \R^d \rightrightarrows \R^d$ be a map such that
   \begin{enumerate}
    \item
     $\phi$ is m-monotone
    \item
     $\abs{y} \geq c \abs{z}$ for some $c > 0$ and all $z$ in some open ball $B_0(\rho)$ around $0$ and $y \in \phi(z)$.
   \end{enumerate}
  Also let $\psi: \R^d \rightrightarrows \R^d$ be any m-monotone map and consider
   \begin{equation}
    A x
     = \mathfrak{A} x,
     \quad
     \dom(\mathfrak{A})
     = \{x \in \dom(\mathfrak{A}): \H x(0) \in \psi(P_1 \H x(0)), \ \H x(1) \in - \phi(P_1 \H x(1)) \}
     \nonumber
   \end{equation}
  where $\mathfrak{A}$ is a first order port-Hamiltonian operator.
  Then $A$ generates an asymptotically stable nonlinear contraction semigroup $(S(t))_{t\geq0}$ on $X = L_2(0,1)^d$ equipped with the inner product $\sp{\cdot}{\cdot}_\H$.
 \end{example}

\section{A Generation Result for Dynamic Controllers}
 \label{sec:generation-dynamic}
 Let us again consider a port-Hamiltonian system $\mathfrak{S} = (\mathfrak{A},\mathfrak{B}, \mathfrak{C})$ of arbitrary order $N \in \N$ which is impedance passive.
 However, instead of a static feedback $\mathfrak{B} x \in - \phi(\mathfrak{C} x)$ we now consider the feedback interconnection $\mathfrak{B} x = - y_c$ and $u_c = \mathfrak{C} x$ with a nonlinear control system $\Sigma_c$.
  \begin{align}
    \left( \begin{array}{c} \frac{\partial}{\partial t} x_c(t) \\ - y_c(t) \end{array} \right)
     &\in M_c \left( \begin{array}{c} x_c(t) \\ u_c(t) \end{array} \right),
     \quad
     t \geq 0.
   \tag{NLC}
  \end{align}
 To motivate the subsequent definitions and assumptions let us first consider the case of a finite dimensional linear system $\Sigma_c = (A_c, B_c, C_c, D_c)$ given by
  \begin{align}
    \frac{\partial}{\partial t} x_c(t) &= A_c x_c(t) + B_c u_c(t)
     \nonumber \\
    y_c(t) &= C_c x_c(t) + D_c u_c(t), \quad t \geq 0.
   \tag{LC}
  \end{align}
 This system is impedance passive if and only if the matrix
  \begin{equation}
   M_c
    := \left( \begin{array}{cc} A_c & B_c \\ - C_c & - D_c \end{array} \right)
    \nonumber
  \end{equation}
 is dissipative (and then m-dissipative since $X_c \times \K^{Nd}$ is finite dimensional and the map is linear).
 So much for the linear and finite dimensional case.
 More general we assume that we have a Hilbert space $X_c$ as the controller state space (equipped with some inner product $\sp{\cdot}{\cdot}_{X_c}$) and $M_c: X_c \times \K^{Nd} \rightrightarrows X_c \times \K^{Nd}$ is a (possibly multi-valued) m-dissipative map for which its minimal section generates a nonlinear s.c. contraction semigroup on $X_c \times \K^{Nd}$.
  \begin{example}
   If we have that $M_c = \left( \begin{smallmatrix} A_c & B_c \\ -C_c & -D_c \end{smallmatrix} \right)$ where $A_c: \dom(A_c) \subseteq X_c \rightrightarrows X_c$ and $-D_c: \dom(-D_c) = \K^{Nd} \rightrightarrows \K^{Nd}$ are m-dissipative and the operators $B_c: \K^{Nd} \rightarrow X_c$ and $C_c: X_c \rightarrow \K^{Nd}$ are assumed to be linear and adjoint to each other, i.e. input and output are collocated, then $M_c: \dom(M_c) = \dom(A_c) \times \K^{Nd} \subseteq X_c \times \K^{Nd} \rightrightarrows X_c \times \K^{Nd}$ is m-dissipative.
   (Note that choosing $B_c = 0$ leads to static feedback as investigated before where the nonlinear system is decoupled from the infinite-dimensional part.)
  \end{example}
 We denote by $\Pi_{X_c}: X_c \times \K^{Nd} \rightarrow X_c$ and $\Pi_{\K^{Nd}}: X_c \times \K^{Nd} \rightarrow \K^{Nd}$ the canonical projections on $X_c$ and $\K^{Nd}$, respectively.
 Then we may define $\A: \dom(\A) \subseteq X \times X_c \rightrightarrows X \times X_c$ as
  \begin{align}
   \A \left( \begin{array}{c} x \\ x_c \end{array} \right)
    &=\left( \begin{array}{c} \mathfrak{A}x \\ \Pi_{X_c} M_c (x_c, \mathfrak{C} x) \end{array} \right)
    \nonumber \\
   \dom(\A)
    &= \{(x,x_c) \in \dom(\mathfrak{A}) \times \Pi_{X_c} \dom(M_c): \ \mathfrak{B} x \in \Pi_{\K^{Nd}} M_c (x_c, \mathfrak{C} x) \}
    \nonumber
  \end{align}
 and we have the following generation theorem for the interconnected system.
 We use the notation $\underline X_c := \overline{\Pi_{X_c} \dom(M_c)}$.
  \begin{theorem}
   \label{thm:generation_nonlinear_controller}
   The (possibly multi-valued) map $\A: \dom(\A) \subset X \times X_c \rightrightarrows X \times X_c$ is m-dissipative on the product space $X \times \underline X_c$,
   thus its minimal section generates a nonlinear s.c. contraction semigroup $(\mathcal{S}(t))_{t\geq0}$ on $X \times \underline X_c$.
  \end{theorem}
 \textbf{Proof.}
 From the Komura-Kato Theorem we know that it is enough to show m-dissipativity.
 Again we may and will assume that $\H = I$.
 First let us show that $\overline{\dom(\A)} = X \times \underline X_c$.
 Take any $(x,x_c) \in X \times \underline X_c$.
 As a first step, let us additionally assume that $x_c \in \Pi_{X_c} \dom(M_c)$.
 Then there are $u_c$ and $y_c \in \K^{Nd}$ such that $(x_c, u_c) \in \dom(M_c)$ and $y_c \in \Pi_{\K^{Nd}} M_c(x_c,u_c)$.
 We need to find a sequence $(x_n)_{n\geq1} \subseteq \dom(\mathfrak{A})$ converging to $x$ (in $X$) and such that $\mathfrak{B}x_n = u_c$ and $\mathfrak{C} x_n = y_c$.
 For this take any $x_0 \in \dom(\mathfrak{A})$ such that $\mathfrak{B} x_0 = y_c$ and $\mathfrak{C} x_0 = u_c$.
 Since $C_c^\infty(0,1)^d$ is dense in $X$ there is a sequence $(z_n)_{n\geq1} \subseteq C_c^\infty(0,1)^d$ converging to $x - x_0$ (in $X$).
 Note that then $x_n := x_0 + z_n \rightarrow x$ with $\mathfrak{B} x_n = \mathfrak{B} x_0 \in \Pi_{\K_{Nd}} M_c(x_c, \mathfrak{C}x_0) = \Pi_{\K_{Nd}} M_c(x_c, \mathfrak{C}x_n)$ does the job.
 As a second step we allow that $(x,x_c)$ merely lies in $X \times \underline X_c$.
 Then there is a sequence $(x_{c,n})_{n\geq1} \subseteq \Pi_{X_c} \dom(M_c)$ such that $\norm{x_{c,n} - x_c}_{X_c} \leq \frac{1}{n}$.
 Further we know from the first step that there are sequences $(x_{n,k}, x_{c,n,k})_{k \geq 1} \subseteq \dom(\A)$ such that $\norm{(x_{n,k},x_{c,n,k}) - (x_n,x_{c,n})}_{X \times X_c} \leq \frac{1}{k}$ and hence the diagonal sequence $(x_{n,n},x_{c,n,n})_{n \geq 1} \subseteq \dom(\A)$ converges to $(x,x_c)$.
 Therefore $\dom(\A)$ is dense in $X \times \underline X_c$.
 Secondly, $\A$ is dissipative since for all $(x,x_c), (\tilde x, \tilde x_c) \in \dom(\A)$, $(\mathfrak{A}x,z_c) \in \A(x,x_c), \ (\mathfrak{A} \tilde x, \tilde z_c) \in \A(\tilde x, \tilde x_c)$ we have
  \begin{align}
   &\Re \sp{(\mathfrak{A}x, z_c) - (\mathfrak{A}\tilde x, \tilde z_c)}{(x,x_c) - (\tilde x, \tilde x_c)}_{X \times X_c}
    \nonumber \\
    &= \Re \sp{\mathfrak{A}(x - \tilde x)}{x - \tilde x}_\H
     + \Re \sp{z_c - \tilde z_c}{x_c - \tilde x_c}_{X_c}
    \nonumber \\
    &\leq \Re \sp{\mathfrak{B}(x - \tilde x)}{\mathfrak{C} (x - \tilde x_c)}_{\K^{Nd}}
     + \Re \sp{z_c - \tilde z_c}{x_c - \tilde x_c}_{X_c}
     \nonumber \\
     &= \Re \sp{\left(\begin{smallmatrix} z_c \\ \mathfrak{B} x \end{smallmatrix} \right) - \left(\begin{smallmatrix} \tilde z_c \\ \mathfrak{B} \tilde x \end{smallmatrix} \right)}{\left(\begin{smallmatrix} x_c \\ \mathfrak{C} x \end{smallmatrix} \right) - \left(\begin{smallmatrix} \tilde x_c \\ \mathfrak{C} \tilde x \end{smallmatrix} \right)}_{X_c \times \K^{Nd}}
     \leq 0
     \nonumber
  \end{align}
 since $\left( \begin{smallmatrix} z_c \\ \mathfrak{B}x \end{smallmatrix} \right) \in M_c \left( \begin{smallmatrix} x_c \\ \mathfrak{C} x \end{smallmatrix} \right)$ and $\left( \begin{smallmatrix} \tilde z_c \\ \mathfrak{B} \tilde x \end{smallmatrix} \right) \in M_c \left( \begin{smallmatrix} \tilde x_c \\ \mathfrak{C} \tilde x \end{smallmatrix} \right)$ for the m-dissipative map $M_c$.
 Finally, we show the range condition $\ran (I -  \A) = X \times X_c$.
 Let $(f,f_c) \in X \times X_c$ be arbitrary.
 We look for $(x,x_c) \in \dom(\A)$ such that
  \begin{equation}
   (x, x_c) - (f,f_c)
    \in \A(x,x_c)
    \nonumber
  \end{equation}
 which may be equivalently expressed as $(x,x_c) \in \dom(\mathfrak{A}) \times \Pi_{X_c} \dom(M_c)$ being such that
  \begin{align}
   (I - \mathfrak{A})x
    &= f
    \nonumber \\
   x_c - f_c
    &\in \Pi_{X_c} M_c \left( \begin{array}{c} x_c \\ \mathfrak{C} x \end{array} \right)
    \nonumber \\
   \mathfrak{B} x
    &\in \Pi_{\K^{Nd}} M_c \left( \begin{array}{c} x_c \\ \mathfrak{C} x \end{array} \right)
    \nonumber
  \end{align}
 where  from the first equality and Lemma \ref{lem:transfer_fct_phs} we get $x = \Phi(1)f + \Psi(1)\mathfrak{B} x$ and $\mathfrak{C} x = F(1)f + G(1) \mathfrak{B} x$.
 Since $G(1)$ is invertible it only remains to solve the problem
  \begin{equation}
   \left( \begin{array}{c} x_c \\ G(1)^{-1} \mathfrak{C} x \end{array} \right) - \left( \begin{array}{c} f_c \\ G(1)^{-1} F(1) f \end{array} \right)
   \in M_c \left( \begin{array}{c} x_c \\ \mathfrak{C}x \end{array} \right).
   \label{eq:generation-thm-dynamic-nl_1}
  \end{equation}
 Since for some $\varepsilon > 0$ small enough $\varepsilon I - \Re G(1)^{-1}$ is still dissipative, clearly  
  \begin{equation}
   \Delta
    := \left( \begin{array}{cc} 0 & \\ & \varepsilon I - G(1)^{-1} \end{array} \right)
    \nonumber
  \end{equation}
 is dissipative and linear from $X_c \times \K^{Nd}$ to $X_c \times \K^{Nd}$ and since $M_c: \dom(M_c) \subseteq X_c \times \K^{Nd} \rightrightarrows X_c \times \K^{Nd}$ is m-dissipative, so is $\Delta + M_c$ by Lemma \ref{lem:perturbation_m-monotone}.
 Hence there is a unique solution $(x_c,\mathfrak{C} x)$ of equation (\ref{eq:generation-thm-dynamic-nl_1}) and we find a unique $(x,x_c) \in \dom(\A)$ such that $(f,f_c) + (x,x_c) \in \A (x,x_c)$.
 \qed

\section{Stabilisation via Nonlinear Dynamic Controllers}
 The idea of this section is to obtain stability results similar to those for the static case, this time in the dynamic controller setup.
 Our results are based on the idea which we employed for the (first) proof of Theorem \ref{thm:exp_stability_nonlinear} where we took $x_0 \in \dom(A)$ and for $x = S(\cdot) x_0$ and some suitable $\eta \in C^1([0,1];\R)$ defined
  \begin{equation}
   \Phi(t)
    = \frac{t}{2} \norm{x(t)}_\H^2 + \sp{x(t)}{\eta P_1^{-1} x(t)}_{L_2}.
    \nonumber
  \end{equation}
 Of course, in the dynamic controller scenario we have to add additional terms corresponding to the finite dimensional controller $\Sigma_c$ as in the preceding section.
 In fact, we assume the following for $\Sigma_c$.
 \begin{assumption}
  \label{asmpt:nonlinear_controller}
   Assume that $\A$ is an m-dissipative operator as in Theorem \ref{thm:generation_nonlinear_controller} and further assume that $0 \in M_c(0)$ and there is $\rho > 0$ and an orthogonal projection $\Pi: \K^{Nd} \rightarrow \K^{Nd}$ on some subspace of $\K^{Nd}$ such that the following hold.
   \begin{enumerate}
    \item
    \begin{align}
     \Re \speq{(z_c,y_c)}{(x_c,u_c)}_{X_c \times \K^{Nd}}
      &\leq - \rho \left( \norm{x_c}_{X_c}^2 + \norm{\Pi u_c}_{\K^{Nd}}^2 \right),
      \nonumber \\
      &\qquad
      (x_c, u_c) \in \dom(M_c), \ (z_c,y_c) \in M_c(x_c,u_c).
      \nonumber 
    \end{align}
    \item
     There is a constant $c > 0$ such that for every $(z_c, -y_c) \in M_c(x_c, u_c)$
      \begin{equation}
       \abs{y_c}^2
        \leq c' \left( \norm{x_c}_{X_c}^2 + \abs{\Pi u_c}^2 \right).
        \label{eq:cond-kernel}
      \end{equation}
    \item
     There are constants $t_0, \delta, c > 0$ such that for every mild solution of the nonlinear control system
      \begin{equation}
       \norm{x_c(t_0)}_{X_c}^2
        \leq (1 - \delta) \norm{x_c(0)}_{X_c}^2
         + c \norm{\Pi u_c}_{L_2(0, t_0)}^2.
         \nonumber
      \end{equation}
   \end{enumerate}
 \end{assumption}
  \begin{remark}
   Note that by time-invariance of the control system the last condition then also holds in the form
      \begin{equation}
       \norm{x_c(t + t_0)}_{X_c}^2
        \leq (1 - \delta) \norm{x_c(t)}_{X_c}^2
         + c \norm{\Pi u_c}_{L_2(t, t + t_0)}^2,
         \quad
         t \geq 0.
         \nonumber
      \end{equation}
   Moreover, the fixed time $t_0 > 0$ may be replaced by $k t_0$ where $k \in \N$ is an arbitrary natural number, i.e. in particular we may choose $t_0 > 0$ as large as we wish.
  \end{remark}
  
  \begin{example}[Collocated case]
  \label{exa:collocated_case}
  One particular case which is covered by the preceding assumption is the following.
  Let $\Sigma_c = (A_c, B_c, C_c, D_c)$ be an impedance passive system with $C_c = B_c^*$ (Hilbert space adjoint with respect to the inner products $\sp{\cdot}{\cdot}_{X_c}$ and $\sp{\cdot}{\cdot}_{\K^{Nd}}$ on $X_c$ and $\K^{Nd}$, respectively) collocated, linear and bounded and $A_c, -D_c$ be m-dissipative.
  Further we assume that
   \begin{enumerate}
    \item
     $A_c(0) = \{0\}$ and there is an equivalent inner product $\speq{\cdot}{\cdot}_{X_c}$ on $X_c$ such that for some $\rho > 0$ and all $x_c \in \dom(A_c), \ z_c \in A_c(x_c)$ one has
      \begin{equation}
       \label{eq:lyapunov}
       \Re \speq{z_c}{x_c}_{X_c}
        \leq - \rho \norm{x_c}_{X_c}^2,
      \end{equation}
    \item
     $0 \in D_c(0), \ \Pi: \K^{Nd} \rightarrow \K^{Nd}$ is an orthogonal projection such that $\abs{w_c} \lesssim \abs{\Pi u_c} \ (u_c \in \dom(D_c), \ w_c \in D_c(u_c))$ and there is $\sigma > 0$ such that for all $x_c \in X_c, \ z_c \in A_c(x_c), \ u_c \in \K^{Nd}$ and $w_c \in D_c(u_c)$ one has
      \begin{equation}
       \label{eq:cond-sip}
       \Re \sp{z_c + B_c u_c}{x_c}_{X_c}
        \leq \Re \sp{C_c x_c + w_c}{u_c}_{\K^{Nd}} - \sigma \abs{w_c}^2.
      \end{equation}
   \end{enumerate}
  Then $M_c = \left( \begin{smallmatrix} A_c & B_c \\ -C_c & -D_c \end{smallmatrix} \right)$ satisfies Assumption \ref{asmpt:nonlinear_controller}.
 \end{example}

 \begin{remark}
  We give some interpretation for the preceding conditions in the collated input/output case.
  The \emph{Lyapunov condition} (\ref{eq:lyapunov}) says that $0$ is a globally exponentially stable equilibrium for the semigroup $(S_c(t))_{t\geq0}$ associated to $A_c$.
  If one has a globally exponentially stable minimum at some other point $x_c^* \in X_c$ one may simply introduce $x_c^{new} := x_c - x_c^*$ as new variable to get to the situation as above.
  (Similar, one may choose a nonzero desired equilibrium $x_c^*$.)
  Conditions (\ref{eq:cond-sip}) and (\ref{eq:cond-kernel}) together may be seen as a strict input passivity condition on the controller system (after getting rid of the redundant parts of the input which only constitute boundary conditions on the system $\mathfrak{S} = (\mathfrak{A}, \mathfrak{B}, \mathfrak{C})$).
  In particular if $D_c = D_c^*$ were linear and symmetric the second condition would read
   \begin{equation}
    \Re \sp{z_c + B_c u_c}{x_c}_{X_c}
     \leq \Re \sp{C_c x_c + D_c u_c}{u_c} - \tilde \sigma \abs{\Pi_{D_c} u_c}^2
     \nonumber
   \end{equation}
  for some $\tilde \sigma > 0$ and $\Pi_{D_c}$ the projection on $\ker{D_c}^\bot$.
 \end{remark}

 We then have the following preliminary, but general result.
 \begin{proposition}
 \label{prop:interconnected_stability}
  Let $\mathfrak{S} = (\mathfrak{A}, \mathfrak{B}, \mathfrak{C})$ be an impedance passive boundary control and observation system and $M_c: \dom(M_c) \subseteq X_c \times \K^{Nd} \rightrightarrows X_c \times \K^{Nd}$ as in Assumption \ref{asmpt:nonlinear_controller}.
  Denote by $(\mathcal{S}(t))_{t\geq0}$ the nonlinear semigroup associated to $\A$ as in Theorem \ref{thm:generation_nonlinear_controller}.
  If there is $q: X \rightarrow \R$ such that $\abs{q(x)} \leq \hat c \norm{x}_\H^2 \ (x \in X)$ and for all solutions $x \in W^1_\infty(\R_+;X) \cap L_\infty(\R_+;\dom(\mathfrak{A}))$ of $\dot x = \mathfrak{A} x$ one has $q(x) \in W_\infty^1(\R_+)$ and
   \begin{equation}
    \norm{x(t)}_\H^2 + \frac{d}{dt} q(x(t))
     \leq c \left( \abs{\mathfrak{B}x(t)}^2 + \abs{\Pi \mathfrak{C} x(t)}^2 \right),
     \quad
     \text{a.e.} \ t \geq 0,
     \nonumber
   \end{equation}
  then $0$ is a globally uniformly exponentially stable equilibrium of $(\mathcal{S}(t))_{t\geq0}$.
 \end{proposition}

%\begin{remark}
%  In the collocated input/output case (see Example \ref{exa:collocated_case}) the first condition means that there is $\tilde c \geq 0$ such that for all $u_c \in \K^{Nd}, w_c \in D_c(u_c)$
%   \begin{equation}
%    \label{eq:cond-kernel}
%    \abs{w_c}
%     \leq \tilde c \norm{B_c u_c}_{X_c}
%   \end{equation}
%  and the last condition may be replaced by
%   \begin{equation}
%    \norm{x(t)}_\H^2 + \frac{d}{dt} q(x(t))
%     \leq c \left( \abs{\mathfrak{B}x(t)}^2 + \abs{D_c^0 \mathfrak{C} x(t)}^2 \right),
%     \quad
%     \text{a.e.} \ t \geq 0,
%     \nonumber
%   \end{equation}
%  where $D_c^0$ denotes the minimal section of $D_c$.
%  To see this let $\Pi: \K^{Nd} \rightarrow \K^{Nd}$ be the projection on $(\ker B_c)^\bot$, then for $u_c \in \dom(D_c)$
%   \begin{equation}
%    \abs{D_c^0 u_c}
%     \lesssim \norm{B_c u_c}_{X_c}
%     \lesssim \abs{\Pi u_c}.
%     \nonumber
%   \end{equation}
% \end{remark}

 \textbf{Proof of Proposition \ref{prop:interconnected_stability}.}
 Since all the maps $S(t): X \times \underline{X_c} \rightarrow X \times \underline{X_c}$ are continuous, it suffices to consider initial data $(x, x_{c,0}) \in \dom(\A)$.
 Moreover, we may and will assume that $t_0 \geq \frac{cc' + c}{2 \sigma} > 0$.
 Let $(x_0, x_{c,0}) \in \dom(\A)$ be arbitrary and let $(x,x_c)(t) := \mathcal{S}(t)(x_0,x_{c,0}) \ (t \geq 0)$ be the unique Lipschitz continuous solution.
 Define
  \begin{equation}
   \Phi(t)
    := t \left( \norm{x(t)}_\H^2 + \norm{x_c(t)}_{X_c}^2 \right) + q(x(t)) + \frac{1 + cc'}{\delta} \int_t^{t+t_0} \norm{x_c}_{X_c}^2 ds,
    \quad
    t \geq 0
    \nonumber
  \end{equation}
 and note that $\frac{d}{dt} (x,x_c)(t) = (\mathfrak{A} x(t), z_c(t)) := \A^0((x,x_c)(t)) \ (\text{a.e.} \ t \geq 0)$ and $\Phi$ is Lipschitz continuous and bounded.
 Then for every $t \geq 2 t_0 > 0$ we have
  \begin{align}
   &\quad
    \Phi(t)
    - \Phi(t_0)
    \nonumber \\
    &= \int_{t_0}^t \frac{d}{ds} \Phi(s) ds
    \nonumber \\
    &= \int_{t_0}^t \norm{(x,x_c)(s)}_{X \times X_c}^2
     + 2 s \Re \sp{\mathfrak{A} x(s)}{x(s)}_X
     + 2 s \Re \sp{z_c(s)}{x_c(s)}_{X_c}
     \nonumber \\
     &\qquad
     + \frac{d}{ds} q(x(s))
     + \frac{1 + cc'}{\delta} (\norm{x_c(s + t_0)}_{X_c}^2 - \norm{x_c(s)}_{X_c}^2) ds
     \nonumber \\
    &\leq \int_{t_0}^t \norm{x_c(s)}_{X_c}^2
     - 2 \sigma s \abs{\Pi \mathfrak{C} x(s)}^2
     + c \left( \abs{\mathfrak{B} x(s)}^2 + \abs{\Pi \mathfrak{C} x(s)}^2 \right)
     \nonumber \\
     &\qquad
     - (1 + cc') \norm{x_c(s)}_{X_c}^2
     + \frac{c(1 + cc')}{\delta} \norm{\Pi \mathfrak{C} x(s)}_{L_2(s,s+t_0)}^2 ds
     \nonumber \\
    &\leq \int_{t_0}^t (- 2 \sigma s + cc' + c) \abs{\Pi \mathfrak{C} x(s)}^2
     \nonumber \\
     &\qquad
     + \frac{c (1 + cc')}{\delta \sigma} \left( \norm{(x,x_c)(s)}_{X \times X_c}^2 - \norm{(x,x_c)(s+t_0)}_{X \times X_c}^2 \right) ds
     \nonumber \\
   &\leq \frac{c (1 + cc')}{\delta \sigma} \left( \norm{(x,x_c)}_{L_2(t_0,2t_0;X \times X_c)}^2 - \norm{(x,x_c)}_{L_2(t,t+t_0;X \times X_c)}^2 \right)
    \nonumber \\
   &\leq t_0 \frac{c (1 + cc')}{\delta \sigma} \norm{(x,x_c)(t_0)}_{X \times X_c}^2.
   \nonumber
  \end{align}
 Since $\frac{\Phi(t)}{t}$ behaves as $\norm{(x,x_c)(t)}^2$ as $t \rightarrow \infty$, we easily deduce exponential stability from this.
 In fact, for $t \geq t_0$ we have
  \begin{align}
   \norm{(x,x_c)(t)}_{X \times X_c}^2
    &= \frac{\Phi(t)}{t} - \frac{q(x(t)) + \norm{x_c(t)}_{X_c}^2}{t} - \frac{1}{t} \norm{x_c}_{L_2(t,t+t_0)}^2
    \nonumber \\
    &\leq \frac{\Phi(t)}{t} + \frac{\hat c \norm{(x,x_c)(t)}_{X \times X_c}}{t}
    \nonumber \\
    &\leq \frac{1}{t} \Phi(t_0) + \frac{\hat c' }{t} \norm{(x,x_c)(t)}_{X \times X_c}^2
    \nonumber \\
    &\leq \frac{t_0}{t}  (1 + \hat c') \norm{(x,x_c)(0)}_{X \times X_c}^2 + \frac{\hat c' }{t} \norm{(x,x_c)(t)}_{X \times X_c}^2,
    \nonumber
  \end{align}
 so that
  \begin{equation}
   \norm{(x,x_c)(t)}_{X \times X_c}^2
    \leq \frac{t_0 (1 + \hat c') \}}{t - \hat c' } \norm{(x,x_c)(0)}_{X \times X_c}^2
    \quad (t > \max \{2 t_0,\hat c' \})
  \end{equation}
 from where exponential stability with constants $M \geq 1$ and $\omega < 0$ independent of $x_0$ follows.
 From density of $\dom(\A)$ in $X \times \underline{X_c}$ and continuity of $\mathcal{S}(t) \ (t \geq 0)$ we conclude
  \begin{equation}
   \norm{\mathcal{S}(t)(x_0,x_{c,0})}_{X \times X_c}
    \leq M e^{\omega t} \norm{(x_0,x_{c,0})}_{X \times X_c},
    \quad (x_0,x_{c,0}) \in X \times \underline{X_c}, \ t \geq 0.
    \nonumber
  \end{equation}
 \qed

\begin{lemma}
\label{lem:q_for_PHS_N=1}
 If $\mathfrak{S} = (\mathfrak{A}, \mathfrak{B}, \mathfrak{C})$ is an impedance passive port-Hamiltonian system of order $N = 1$, then there is $q \in C^1(X;\R)$ with $\abs{q(x)} \leq \hat c \norm{x}_\H^2 \ (x \in X)$ such that for every solution $x \in W_\infty^1(\R_+;X) \cap L_\infty(\R_+;\dom(\mathfrak{A}))$ of $\dot x = \mathfrak{A} x$ one has
  \begin{equation}
   \norm{x(t)}_\H^2 + \frac{d}{dt} q(x(t))
    \leq c \abs{(\H x(t))(1)}^2,
     \quad
     \text{a.e} \ t \geq 0.
     \nonumber
  \end{equation}
\end{lemma}
 \textbf{Proof.}
 Take $q(x) = \sp{x}{\eta P_1^{-1} x}_{L_2}$ as in the proof of Theorem \ref{thm:exp_stability_nonlinear} and the calculations made there confirm the assertion.
 \qed

Thus Proposition \ref{prop:interconnected_stability} and Lemma \ref{lem:q_for_PHS_N=1} together say the following.

\begin{theorem}
 \label{thm:exp_stability_N=1}
 Let $\mathfrak{S} = (\mathfrak{A}, \mathfrak{B}, \mathfrak{C})$ be an impedance passive port-Hamiltonian system of order $N = 1$ and $M_c: \dom(M_c) \subseteq X_c \times \K^{Nd} \rightrightarrows X_c \times \K^{Nd}$ be as in Assumption \ref{asmpt:nonlinear_controller}.
 Further assume that
  \begin{equation}
   \abs{(\H x)(1)}^2
    \lesssim \abs{\mathfrak{B}x}^2 + \abs{\Pi \mathfrak{C}x}^2,
    \quad x \in \dom(\mathfrak{A}).
    \nonumber
  \end{equation}
 Then the interconnected map $\A$ from Theorem \ref{thm:generation_nonlinear_controller} generates a s.c. contraction semigroup $(S(t))_{t\geq0}$ on $X \times \underline X_c$ with globally exponentially stable equilibrium 0.
\end{theorem}

\section{Stabilisation of Second Order Systems}
 In this section we aim for a generalisation of Theorem \ref{thm:exp_stability_N=1} to the case where
  \begin{equation}
   \mathfrak{A} x
    = P_2 (\H x)'' + P_1 (\H x)' + P_0 (\H x)
    \nonumber
  \end{equation}
 is a port-Hamiltonian operator of second order ($N = 2$).
 Again we assume that $\H \in W_\infty^1(0,1)^{d \times d}$ is Lipschitz continuous.
 For the case of (static and dynamic) linear feedback stabilisation, see \cite{AugnerJacob_2014}, where exponential stability has been proved under the assumption
  \begin{equation}
   \abs{(\H x)(0)}^2 + \abs{(\H x)(1)}^2 + \abs{(\H x)'(0)}^2
    \lesssim \abs{\mathfrak{B}x}^2 + \abs{\Pi \mathfrak{C}x}^2,
    \quad
    x \in \dom(\mathfrak{A}).
    \nonumber
  \end{equation}
 Of course, the proof there used the Gearhart-Greiner-Pr\"uss Theorem, so lacks any possible generalisation to the nonlinear scenario.
 Therefore we aim to apply Proposition \ref{prop:interconnected_stability} which amounts to finding a suitable $q \in C^1(X;\R)$ satisfying the assumptions of Proposition \ref{prop:interconnected_stability}.

 Unfortunately, it is very hard to prove existence of such a functional $q$ without any further restrictions on $\H$ and the matrices $P_0$ and $P_1$ and, in fact, we did not succeed in proving the general result, but had to impose further constraints on the matrix-valued function $\H$ and the matrices $P_0$ and $P_1$.

 \begin{lemma}
  \label{lem:N=2}
  Let $\mathfrak{S} = (\mathfrak{A}, \mathfrak{B}, \mathfrak{C})$ be an impedance passive port-Hamiltonian system of order $N = 2$.
  Further assume that $\H', P_0$ and $P_1$ are small compared to $\H$, i.e.
    \begin{align}
     2
      &> \norm{(\H' \H^{-1} + P_2^{-1} P_1)(\z - 1)}_{L_\infty(0,1;\K^{d \times d})}
       \nonumber \\
       &\quad
       + \frac{1}{\sqrt{2}} \norm{(P_0^* P_2^{-1} + P_2^{-1} P_0 - P_1 P_2^{-1} P_1 P_2^{-1})}
       \nonumber \\
       &\quad
       + \frac{1}{2} \norm{((P_2^{-1} P_0)^* P_2^{-1} P_0)}.
    \end{align}
  Then there is $q: X \rightarrow \R$ such that $\abs{q(x)} \leq \hat c \norm{x}_\H^2 \ (x \in X)$ and for all solutions $x \in W_\infty^1(\R_+;X) \cap L_\infty(\R_+;\dom(\mathfrak{A}))$ of $\dot x = \mathfrak{A} x$ the function $q(x)$ lies in $W^1_{\infty,loc}(\R_+)$ with
   \begin{equation}
    \norm{x(t)}_\H^2 + \frac{d}{dt} q(x(t))
     \leq c \left( \abs{(\H x)(t,0)}^2 + \abs{(\H x)'(t,0)}^2 + \abs{(\H x)(t,1)}^2 \right)
     \nonumber
   \end{equation}
  for a.e. $t \geq 0$.
 \end{lemma}
 \textbf{Proof.}
 We define the real-valued functional $q: X \rightarrow \R$ as
  \begin{align}
   q(x)
    &:= \Re \sp{x}{\eta P_2^{-1} \int_0^\cdot x(\xi) d\xi}_{L_2}
    \nonumber \\
    &\quad
    - \frac{1}{2} \sp{P_2^{-1} \int_0^\cdot x(\xi) d\xi}{\eta P_1 P_2^{-1} \int_0^\cdot x(\xi) d\xi}_{L_2},
    \quad
    x \in X
  \end{align}
 where the scalar function $\eta \in C^\infty([0,1];\R)$ may be chosen suitable at a later point.
 Then for every Lipschitz continuous solution $x \in W_\infty^1(\R_+;X) \cap L_\infty(\R_+;\dom(\mathfrak{A}))$ of the evolution equation $\dot x = \mathfrak{A} x$ we obtain (omitting the parameter $t$ for brevity and employing Lemma \ref{lem:real_part}) that
  \begin{align}
   &\quad
    \frac{d}{dt} q(x)
    \nonumber \\
    &= \Re \sp{P_2 (\H x)''}{\eta P_2^{-1} \int_0^\cdot x(\xi) d\xi}_{L_2}
     + \Re \sp{x}{\eta \int_0^\cdot (\H x)''(\xi) d\xi}_{L_2}
     \nonumber \\
     &\quad
     + \Re \sp{P_1 (\H x)' + P_0 (\H x)}{\eta P_2^{-1} \int_0^\cdot x(\xi) d\xi}_{L_2}
     \nonumber \\
     &\quad
     + \Re \sp{x}{\eta P_2^{-1} \int_0^\cdot P_1 (\H x)'(\xi) + P_0 (\H  x)(\xi) d\xi}_{L_2}
     \nonumber \\
     &\quad
     - \Re \sp{\int_0^\cdot (\H x)''(\xi) d\xi}{\eta P_1 P_2^{-1} \int_0^\cdot x(\xi) d\xi}_{L_2}
     \nonumber \\
     &\quad
     - \Re \sp{P_2^{-1} \int_0^\cdot P_1 (\H x)'(\xi) + P_0 (\H x)(\xi) d\xi}{\eta P_1 P_2^{-1} \int_0^\cdot x(\xi) d\xi}_{L_2}
     \nonumber \\
    &= 2 \Re \sp{(\H x)'}{\eta x}_{L_2}
     + \Re \sp{(\H x)'}{\eta' \int_0^\cdot x(\xi) d\xi}_{L_2}
     \nonumber \\
     &\quad
     - \Re \eta(1) \sp{(\H x)'(1)}{\int_0^1 x(\xi) d\xi}_{\K^d}
     - \Re \sp{(\H x)'(0)}{\eta x}_{L_2}
     \nonumber \\
     &\quad
     + \Re \sp{P_1 (\H x)'}{\eta P_2^{-1} \int_0^\cdot x(\xi) d\xi}_{L_2}
     + \Re \sp{x}{\eta \int_0^\cdot P_2^{-1} P_1 (\H x)'(\xi) d\xi}_{L_2}
     \nonumber \\
     &\quad
     - \Re \sp{\int_0^\cdot x(\xi) d\xi}{\eta P_2^{-1} P_0 (\H x)}_{L_2}
     + \Re \sp{x}{\eta P_2^{-1} P_0 \int_0^\cdot (\H x)(\xi) d\xi}_{L_2}
     \nonumber \\
     &\quad
     - \Re \sp{P_1 (\H x)'}{\eta P_2^{-1} \int_0^\cdot x(\xi) d\xi}_{L_2}
     + \Re \sp{P_1 (\H x)'(0)}{\eta P_2^{-1} \int_0^\cdot x(\xi) d\xi}_{L_2}
     \nonumber \\
     &\quad
     - \Re \sp{P_2^{-1} P_1 (\H x)}{\eta P_1 P_2^{-1} \int_0^\cdot x(\xi) d\xi}_{L_2}
     \nonumber \\
     &\quad
     + \Re \sp{P_2^{-1} P_1 (\H x)(0)}{\eta P_1 P_2^{-1} \int_0^\cdot x(\xi) d\xi}_{L_2}
     \nonumber \\
     &\quad
     - \Re \sp{\int_0^\cdot P_2^{-1} P_0 (\H x)(\xi)}{\eta P_1 P_2^{-1} \int_0^\cdot x(\xi) d\xi}_{L_2}
     \nonumber \\
    &= \sp{(\eta \H' - 2 \eta' \H) x}{x}_{L_2}
     + \left[ \sp{x(\z)}{(\eta \H)(\z) x(\z)}_{\K^d} \right]_0^1
     \nonumber \\
     &\quad
     - \Re \sp{\H x}{\eta '' \int_0^\cdot x(\xi) d\xi}_{L_2}
     + \eta'(1) \Re \sp{(\H x)(1)}{\int_0^1 x(\xi) d\xi}_{\K^d}
     \nonumber \\
     &\quad
     - \eta(1) \Re \sp{(\H x)'(1)}{\int_0^1 x(\xi) d\xi}_{\K^d}
     - \Re \sp{x}{\eta (\H x)'(0)}_{L_2}
     \nonumber \\
     &\quad
     + \Re \sp{x}{\eta P_2^{-1} P_1 (\H x)}_{L_2}
     - \Re \sp{x}{\eta P_2^{-1} P_1 (\H x)(0)}_{L_2}
     \nonumber \\
     &\quad
     - \Re \sp{\int_0^\cdot x(\xi) d\xi}{\eta P_2^{-1} P_0 (\H x)}_{L_2}
     + \Re \sp{x}{\eta P_2^{-1} P_0 \int_0^\cdot (\H x)(\xi) d\xi}_{L_2}
     \nonumber \\
     &\quad
     + \Re \sp{P_1 (\H x)'(0)}{\eta P_2^{-1} \int_0^\cdot x(\xi) d\xi}_{L_2}
     \nonumber \\
     &\quad
     - \Re \sp{P_2^{-1} P_1 (\H x)}{\eta P_1 P_2^{-1} \int_0^\cdot x(\xi) d\xi}_{L_2}
     \nonumber \\
     &\quad
     + \Re \sp{P_2^{-1} P_1 (\H x)(0)}{\eta P_1 P_2^{-1} \int_0^\cdot x(\xi) d\xi}_{L_2}
     \nonumber \\
     &\quad
     - \Re \sp{\int_0^\cdot P_2^{-1} P_0 (\H x)(\xi)}{\eta P_1 P_2^{-1} \int_0^\cdot x(\xi) d\xi}_{L_2}
     \nonumber \\
    &\leq
     \sp{(\varepsilon \H + \eta \H' - 2 \eta' \H + \eta \Re (P_2^{-1} P_1 \H)) x}{x}_{L_2}
     \nonumber \\
     &\quad
      + \Re \sp{\H x}{(- \eta'' + P_0^* P_2^{-1} \eta + P_2^{-1} P_0 \eta - P_1 P_2^{-1} P_1 P_2^{-1} \eta) \int_0^\cdot x(\xi) d\xi}_{L_2}
      \nonumber \\
      &\quad
      - \Re \sp{P_2^{-1} P_0 \int_0^\cdot \H x(\xi)}{\eta P_2^{-1} P_0 \int_0^\cdot x(\xi) d\xi}_{L_2}
      \nonumber \\
      &\quad
      + c_{\varepsilon, \eta} \left( (1 + \abs{\eta(0)}^2) \abs{(\H x)(0)}^2
       + \abs{(\H x)'(0)}^2
       \right.
       \nonumber \\
       &\qquad
       \left.
       + \abs{\eta(1)}^2 \abs{(\H x)(1)}^2
       + \abs{\eta'(1)}^2 \abs{(\H x)(1)}^2 \right)
  \end{align}
 for every $\varepsilon > 0$ and a constant $c_{\varepsilon, \eta} > 0$ which may depend on $\varepsilon > 0$ and $\eta$, but which is independent of $x$.
 We now estimate in the following ways.
 On the one hand
  \begin{align}
   &\quad
   \Re \sp{\H x}{(- \eta'' + P_0^* P_2^{-1} \eta + P_2^{-1} P_0 \eta + P_2^{-1} P_0 \eta - P_1 P_2^{-1} P_1 PO_2^{-1} \eta) \int_0^\cdot x(\xi) d\xi}_{L_2}
   \nonumber \\
   &\leq \norm{\H x}_{L_2} \norm{- \eta'' + (P_0^* P_2^{-1} + P_2^{-1} P_0 - P_1 P_2^{-1} P_1 P_2^{-1}) \eta}_{L_\infty(0,1;\K^{d \times d})} \norm{\int_0^\cdot x(\xi) d\xi}_{L_2}
   \nonumber \\
   &\leq \norm{\H x}_{L_2} \norm{- \eta'' + (P_0^* P_2^{-1} + P_2^{-1} P_0 - P_1 P_2^{-1} P_1 P_2^{-1}) \eta}_{L_\infty(0,1;\K^{d \times d})} \frac{1}{\sqrt{2}} \norm{x}_{L_2}
   \nonumber
  \end{align}
 and on the other hand
  \begin{align}
   &\quad
    - \Re \sp{P_2^{-1} P_0 \int_0^\cdot (\H x)(\xi) d\xi}{\eta P_2^{-1} P_0 \int_0^\cdot x(\xi) d\xi}_{L_2}
    \nonumber \\
    &\leq \norm{(P_2^{-1} P_0)^* P_2^{-1} P_0 \eta}_{L_\infty(0,1;\K^{d \times d})} \norm{\int_0^\cdot (\H x)(\xi) d\xi}_{L_2} \norm{\int_0^\cdot x(\xi) d\xi}_{L_2}
    \nonumber \\
    &\leq \norm{(P_2^{-1} P_0)^* P_2^{-1} P_0 \eta}_{L_\infty(0,1;\K^{d \times d})} \frac{1}{2} \norm{\H x}_{L_2} \norm{x}_{L_2}.
    \nonumber
  \end{align}
 Therefore,
  \begin{align}
   &\quad
    \frac{d}{dt} q(x)
    \nonumber \\
    &\leq \sp{\left[ \varepsilon - 2 \eta' + \eta (\H' \H^{-1} + 2 \Re (P_2^{-1} P_1 \H) \H^{-1}) \right] \H x}{x}_{L_2}
     \nonumber \\
     &\quad
     - \left[ \frac{\norm{- \eta'' + (P_0^* P_2^{-1} + P_2^{-1} - P_1 P_2^{-1} P_1 P_2^{-1}) \eta}_{L_\infty(0,1;\K^{d \times d})}}{\sqrt{2}}
     \right.
     \nonumber \\
     &\quad
     \left. + \frac{\norm{(P_2^{-1} P_0)^* P_2^{-1} P_0 \eta}_{L_\infty(0,1;\K^{d \times d})}}{2} \right] \norm{\H x}_{L_2} \norm{x}_{L_2}
     \nonumber
  \end{align}
 and we only have to find a suitable function $\eta$ such that
  \begin{align}
   2 \eta'
    &\geq \varepsilon
     + \norm{(\H' \H^{-1} + P_2^{-1} P_1) \eta}_{L_\infty(0,1;\K^{d \times d})}
     \nonumber \\
     &\quad
     + \frac{1}{\sqrt{2}} \norm{- \eta'' + (P_0^* P_2^{-1} + P_2^{-1} P_0 - P_1 P_2^{-1} P_1 P_2^{-1}) \eta}_{L_\infty(0,1;\K^{d \times d})}
     \nonumber \\
     &\quad
     + \frac{1}{2} \norm{(P_2^{-1} P_0)^* P_2^{-1} P_0 \eta}_{L_\infty(0,1;\K^{d \times d})}.
     \nonumber
  \end{align}
 In particular for the choice $\eta(\z) = 1 - \z$ we obtain the condition
  \begin{align}
   &\norm{(\H' \H^{-1} + P_2^{-1} P_1)(\z - 1)}_{L_\infty(0,1;\K^{d \times d})}
    \nonumber \\
    &+ \frac{1}{\sqrt{2}} \norm{P_0^* P_2^{-1} + P_2^{-1} P_0 - P_1 P_2^{-1} P_1 P_2^{-1}}
    \nonumber \\
    &+ \frac{1}{2} \norm{(P_2^{-1} P_0)^* P_2^{-1} P_0}
    \nonumber \\
    &\quad
    < 2.
    \nonumber
  \end{align}
 The assertion follows.
 \qed

 The interplay of Lemma \ref{lem:N=2} with Proposition \ref{prop:interconnected_stability} then implies the following.

\begin{theorem}
 \label{thm:exp_stability_N=2}
 Let $\mathfrak{S} = (\mathfrak{A}, \mathfrak{B}, \mathfrak{C})$ be an impedance passive port-Hamiltonian system of order $N = 2$ and $M_c: \dom(M_c) \subseteq X_c \times \K^{Nd} \rightrightarrows X_c \times \K^{Nd}$ as in Assumption \ref{asmpt:nonlinear_controller}.
 Further assume that
  \begin{align}
     2
      &> \norm{(\H' \H^{-1} + P_2^{-1} P_1)(\z - 1)}_{L_\infty(0,1;\K^{d \times d})}
       \nonumber \\
       &\quad
       + \frac{1}{\sqrt{2}} \norm{(P_0^* P_2^{-1} + P_2^{-1} P_0 - P_1 P_2^{-1} P_1 P_2^{-1})}
       \nonumber \\
       &\quad
       + \frac{1}{2} \norm{((P_2^{-1} P_0)^* P_2^{-1} P_0)}.
  \end{align}
 and
  \begin{equation}
   \abs{(\H x)(0)}^2 + \abs{(\H x)(1)}^2 + \abs{(\H x)'(0)}^2
    \lesssim \abs{\mathfrak{B}x}^2 + \abs{\Pi \mathfrak{C}x}^2,
    \quad x \in \dom(\mathfrak{A}).
    \nonumber
  \end{equation}
 Then the interconnected map $\A$ from Theorem \ref{thm:generation_nonlinear_controller} generates a s.c. contraction semigroup $(S(t))_{t\geq0}$ on $X \times \underline X_c$ with globally exponentially stable equilibrium 0.
\end{theorem}

\section{Stabilisation of the Euler-Bernoulli Beam}
 We investigate how the general result Proposition \ref{prop:interconnected_stability} may be used to design exponentially stabilising controllers for the Euler-Bernoulli Beam equation, i.e. the dynamical system
  \begin{equation}
   \rho(\z) \omega_{tt}(t,\z) + (EI \omega_{\z\z})_{\z\z}(t,\z),
    \quad \z \in (0,1), \ t \geq 0.
    \label{eqn:eb}
  \end{equation}
 The energy of the system is given by
  \begin{equation}
   E(t)
    = \frac{1}{2} \int_0^1 \rho(\z) \abs{\omega_t(t,\z)}^2 + EI(\z) \abs{\omega_{\z\z}(t,\z)}^2 d\z,
    \quad t \geq 0.
    \nonumber
  \end{equation}
 To rewrite (\ref{eqn:eb}) as a (second-order) port-Hamiltonian system, we set
  \begin{align}
   x
    &= \left( \begin{array}{c} x_1 \\ x_2 \end{array} \right)
    := \left( \begin{array}{c} \rho \omega_t \\ \omega_{\z\z} \end{array} \right),
    \qquad
   \H
    = \left( \begin{array}{cc} \H_1 & \\ & \H_2 \end{array} \right)
    := \left( \begin{array}{cc} \rho^{-1} & \\ & EI \end{array} \right),
    \nonumber \\
   P_2
    &= \left( \begin{array}{cc} & -P^* \\ P & \end{array} \right)
    := \left( \begin{array}{cc} & -1 \\ 1 & \end{array} \right)
    \nonumber
  \end{align}
 and $P_0 = P_1 = 0$, so (\ref{eqn:eb}) takes the form
  \begin{equation}
   \frac{\partial}{\partial t} x(t,\z)
    = P_2 \frac{\partial^2}{\partial \z^2} \H x(t,\z)
    =: \mathfrak{A} x(t,\z),
    \quad \z \in (0,1), \ t \geq 0.
    \nonumber
  \end{equation}
 Note that
  \begin{align}
   \Re \sp{\mathfrak{A}x}{x}_\H
    &= \Re \left[ (\H_1 x_1)'(\z)^* (\H_2 x_2)(\z) - (\H_1 x_1)(\z)^* (\H_2 x_2)'(\z) \right]_0^1
    \nonumber \\
    &= \Re \left[ \omega_{t\z}(\z)^* (EI \omega_{\z\z})(\z) - \omega_t(\z)^* (EI \omega_{\z\z})_\z (\z) \right]_0^1
    \nonumber
  \end{align}
 for all $x = \left( \begin{smallmatrix} \rho \omega_t \\ \omega_{\z\z} \end{smallmatrix} \right) \in \dom(\mathfrak{A})$.
 From here, several choices of $\mathfrak{B}$ and $\mathfrak{C}$ are possible to make $\mathfrak{S} = (\mathfrak{A}, \mathfrak{B}, \mathfrak{C})$ an impedance passive port-Hamiltonian boundary control and observation system.
 In that case (and if $\rho^{-1}, EI \in L_\infty(0,1)$ are uniformly positive) for any $m$-monotone $\phi: \K^2 \rightrightarrows \K^2$ the operator $A = \mathfrak{A}|_{\dom(A)}, \ \dom(A) = \{ x \in \dom(\mathfrak{A}): \mathfrak{B} x \in - \phi(\mathfrak{C} x) \}$ generates a s.c. contraction semigroup on $X = L_2(0,1)^2$ (which is a $C_0$-semigroup if $\phi \in \K^{2 \times 2}$ is linear).
 Lemma \ref{lem:N=2} gives some conditions under which the system can be exponentially stabilised, however these conditions are rather strong and the proof of Lemma \ref{lem:N=2} does not take into account the additional structure of the Euler-Bernoulli beam, in particular those of the matrices $P_i$.
 We therefore give a result analogous to Lemma \ref{lem:N=2} making use of the Euler-Bernoulli beam structure.
  \begin{lemma}
   \label{lem:eb}
   Assume that $\mathfrak{S} = (\mathfrak{A}, \mathfrak{B}, \mathfrak{C})$ is an impedance passive second order port-Hamiltonian system of the form
    \begin{equation}
     \H
      = \left( \begin{array}{cc} \H_1 & \\ & \H_2 \end{array} \right),
      \quad
     P_2
      = \left( \begin{array}{cc} & -P^* \\ P & \end{array} \right),
      \quad
      P_1 = P_0 = 0.
      \nonumber
    \end{equation}
   Further assume that $\H_i \in W_\infty^1(0,1;\K^{d/2 \times d/2}) \ (i = 1, 2)$ where $d \in 2 \N$ is even and
    \begin{align}
     \sup \left\{ \norm{\H_2' \H_2^{-1}}_{L_\infty},
      \norm{\H_1' \H_1^{-1}}_{L_\infty} \right\}
       < 1.
       \nonumber
    \end{align}
   Then there is $q \in C^1(X;\R)$ with $\abs{q(x)} \leq \hat c \norm{x}_\H^2 \ (x \in X)$ such that for all solutions $x \in W_\infty^1(\R_+;X) \cap L_\infty(\R_+;\dom(\mathfrak{A}))$ of $\dot x = \mathfrak{A} x$ one has
    \begin{equation}
     \norm{x(t)}_\H^2 + \frac{d}{dt} q(x(t))
      \leq c \left( \abs{(\H x)(0)}^2 + \abs{(\H_1 x_1)'(0)}^2 + \abs{(\H_2 x_2)(1)}^2 \right),
      \quad \text{a.e.} \ t \geq 0.
      \nonumber
    \end{equation}
  \end{lemma}
 
% \begin{remark}
%  By symmetry, one also gets the following estimates in Lemma \ref{lem:eb}.
%   \begin{align}
%    \norm{x(t)}_\H^2 + \frac{d}{dt} q(x(t))
%     &\leq c \left( \abs{(\H x)(1)}^2 + \abs{(\H_1 x_1)'(1)}^2 + \abs{(\H_2 x_2)(0)}^2 \right),
%     \nonumber \\
%    \norm{x(t)}_\H^2 + \frac{d}{dt} q(x(t))
%     &\leq c \left( \abs{(\H x)(0)}^2 + \abs{(\H_2 x_2)'(0)}^2 + \abs{(\H_1 x_1)(1)}^2 \right),
%     \nonumber \\
%    \norm{x(t)}_\H^2 + \frac{d}{dt} q(x(t))
%     &\leq c \left( \abs{(\H x)(1)}^2 + \abs{(\H_2 x_2)'(1)}^2 + \abs{(\H_1 x_1)(0)}^2 \right).
%     \nonumber
%   \end{align}
% \end{remark}
  \textbf{Proof.}
  This proof is based on the technique used in \cite{ChenEtAl_1987} for a chain of Euler-Bernoulli beams with the particular boundary condition $\omega_t(0) = \omega_\z(0) = 0$ at the left end.
  Let $\eta \in C^2([0,1];\R)$ be a twice continuously differentiable real function which we choose at a later point and define
   \begin{equation}
    q(x)
     := \Re \sp{x_1}{\eta P^{-1} \int_0^\cdot x_2(\xi) d\xi},
     \quad x = (x_1, x_2) \in X.
     \nonumber
   \end{equation}
  Then for every Lipschitz continuous solution $x \in W_\infty^1(\R_+;X) \cap L_\infty(\R_+;\dom(\mathfrak{A}))$ of the evolution equation $\dot x = \mathfrak{A} x$ we have (using Lemma \ref{lem:real_part} again)
   \begin{align}
    &\frac{d}{dt} q(x)
     \nonumber \\
     &= \Re \sp{P^{-*} x_{1,t}}{\eta \int_0^\cdot x_2(\xi) d\xi}_{L_2}
      + \Re \sp{x_1}{\eta \int_0^\cdot P^{-1} x_{2,t}(\xi) d\xi}_{L_2}
      \nonumber \\
     &= - \Re \sp{(\H_2 x_2)''}{\eta \int_0^\cdot x_2(\xi) d\xi}_{L_2}
      + \Re \sp{x_1}{\eta \int_0^\cdot (\H_1 x_1)''(\xi) d\xi}_{L_2}
      \nonumber \\
     &= \Re \sp{(\H_2 x_2)'}{\eta x_2}_{L_2}
      + \Re \sp{(\H_2 x_2)'}{\eta' \int_0^\cdot x_2(\xi) d\xi}_{L_2}
      \nonumber \\
      &\quad
      - \Re \sp{\eta(1) (\H_2 x_2)'(1)}{\int_0^1 x_2(\xi) d\xi}_{\K^{d/2}}
      \nonumber \\
      &\quad
      + \Re \sp{x_1}{\eta (\H_1 x_1)'}_{L_2}
      - \Re \sp{\eta x_1}{(\H_1 x_1)'(0)}_{L_2}
      \nonumber \\
     &= - \frac{1}{2} \sp{x_2}{((\eta \H_2)'-2 \eta \H_2') x_2}_{L_2}
      + \frac{1}{2} \left[ \sp{x_2(\z)}{(\eta \H_2)(\z) x_2(\z)}_{\K^{d/2}} \right]_0^1
      \nonumber \\
      &\quad
      - \Re \sp{\H_2 x_2}{\eta'' \int_0^\cdot x_2(\xi) d\xi}_{L_2}
      - \Re \sp{\H_2 x_2}{\eta' x_2}_{L_2}
      \nonumber \\
      &\quad
      + \Re \sp{\eta'(1) (\H_2 x_2)(1)}{\int_0^1 x_2(\xi) d\xi}_{\K^{d/2}}
      - \Re \sp{\eta(1) (\H_2 x_2)'(1)}{\int_0^1 x_2(\xi) d\xi}_{\K^{d/2}}
      \nonumber \\
      &\quad
      - \frac{1}{2} \sp{x_1}{((\eta \H_1)'- 2 \eta \H_1') x_1}_{L_2}
      + \frac{1}{2} \left[ \sp{x_1(\z)}{(\eta \H_1)(\z) x_1(\z)}_{\K^{d/2}} \right]_0^1
      \nonumber \\
      &\quad
      - \Re \sp{\eta x_1}{(\H_1 x_1)'(0)}_{L_2}
      \nonumber \\
     &\leq \frac{1}{2} \sp{(- \eta' - \eta \H_2' \H_2^{-1} + \varepsilon) \H_2 x_2}{x_2}_{L_2}
      \nonumber \\
      &\quad
      + \frac{1}{2} \sp{(- \eta' + \eta \H_1' \H_1^{-1} + \varepsilon) \H_1 x_1}{x_1}_{L_2}
      - \Re \sp{\H_2 x_2}{\eta'' \int_0^\cdot x_2(\xi) d\xi}_{L_2}
      \nonumber \\
      &\quad
      + c_{\varepsilon}
       \left( \abs{\eta(1) (\H_2 x_2)'(1)
        - \eta'(1) (\H_2 x_2)(1)}^2
        + \norm{\eta}_{L_\infty} \abs{(\H_1 x_1)'(0)}^2
        \right.
        \nonumber \\
        &\qquad
        \left.
        + \abs{\eta(0)} \abs{(\H_1 x_1)(0)}^2
        + \abs{\eta(1)} \abs{(\H_1 x_1)(1)}^2 \right).
      \nonumber
   \end{align}
  We therefore need to find $\eta$ such that $\eta(1) = 0$ and the following conditions hold true for some $\varepsilon > 0$.
   \begin{align}
    \eta'
     &\leq - \left( \abs{\eta} \norm{\H_2' \H_2^{-1}}_{L_\infty} + \norm{\eta'' \sqrt{\z}}_{L_\infty} + \varepsilon \right)
     \nonumber \\
    \eta'
     &\leq - \left( \abs{\eta} \norm{\H_1' \H_1^{-1}}_{L_\infty} + \varepsilon \right).
     \nonumber
   \end{align}
  This gives the assertion of the lemma under the condition on $\H$, using the simple choice $\eta(\z) = 1 - \z$.
  \qed

\begin{theorem}
 \label{thm:exp_stability_eb}
 Let $\mathfrak{S} = (\mathfrak{A}, \mathfrak{B}, \mathfrak{C})$ be an impedance passive port-Hamiltonian system of order $N = 2$ of Euler-Bernoulli type as in Lemma \ref{lem:eb} and $M_c: \dom(M_c) \subseteq X_c \times \K^{Nd} \rightrightarrows X_c \times \K^{Nd}$ as in Assumption \ref{asmpt:nonlinear_controller}.
 Further assume that
  \begin{equation}
   \abs{(\H x)(0)}^2 + \abs{(\H_1 x_1)'(0)}^2 + \abs{(\H_2 x_2)(1)}^2
    \lesssim \abs{\mathfrak{B}x}^2 + \abs{\Pi \mathfrak{C}x}^2,
    \quad x \in \dom(\mathfrak{A}).
    \nonumber
  \end{equation}
 Then the interconnected map $\A$ from Theorem \ref{thm:generation_nonlinear_controller} generates a s.c. contraction semigroup $(S(t))_{t\geq0}$ on $X \times \underline X_c$ with globally uniformly exponential stable equilibrium 0.
\end{theorem}
 \textbf{Proof.}
 Combine Lemma \ref{lem:eb} with Proposition \ref{prop:interconnected_stability}.
 \qed

\section*{Acknowledgements}
The author would like to thank Birgit Jacob for arousing his interest in this research topic, for fruitful discussions and helpful advice as well as for careful proof-reading.
He also is grateful to Hans Zwart for making a research visit at the University of Twente possible and for several discussions leading to a better understanding of the topic and also pointing out a mistake in an earlier version of this article.

\end{document}